\documentclass[12pt,a4paper]{article}%
\usepackage{amsmath}
\usepackage{graphicx}
\usepackage{a4}%
\usepackage{amsfonts}%
\usepackage{amssymb}

\begin{document}

\title{Algorithm for multiplying Schubert classes}
\author{Haibao Duan and Xuezhi Zhao\\Institute of Mathematics, Chinese Academy of Sciences, \\Beijing 100080, dhb@math.ac.cn\\Department of Mathematics, Capital Normal University\\Beijing 100037, zhaoxve@mail.cnu.edu.cn}
\date{\ \ }
\maketitle

\begin{abstract}
Based on the multiplicative rule of Schubert classes obtained in [Du$_{3}$],
we present an algorithm computing the product of two arbitrary Schubert
classes in a flag variety $G/H$, where $G$ is a compact connected Lie group
and $H\subset G$ is the centralizer of a one-parameter subgroup in $G$.

Since all Schubert classes on $G/H$ constitute an basis for the integral
cohomology $H^{\ast}(G/H)$, the algorithm gives a method to compute the
cohomology ring $H^{\ast}(G/H)$ independent of the classical spectral sequence
method due to Leray [L$_{1}$,L$_{2}$] and Borel [Bo$_{1}$, Bo$_{2}$].

\begin{description}
\item 2000 Mathematical Subject Classification: 14N15; 14M10 (55N33; 22E60).

\item Key words and phrases: flag manifolds; Schubert varieties; cohomology;
Cartan matrix.

\end{description}
\end{abstract}

\section{\textbf{Introduction}\ \ \ }

This paper presents an algorithm computing the integral cohomology ring of a
flag manifold $G/H$, where $G$ is a compact connected Lie group and $H\subset
G$ is the centralizer of a one-parameter subgroup.

The determination of the integral cohomology of a topological space is a
classical problem in algebraic topology. However, since a flag manifold $G/H$
is canonically an algebraic variety whose Chow ring is isomorphic to the
integral cohomology $H^{\ast}(G/H)$, a complete description for the ring
$H^{\ast}(G/H)$ is also of fundamental importance to the algebraic
intersection theory of $G/H$ ([K,S$_{2}$]).

\bigskip

In general, an entire account for the integral cohomology $H^{\ast}(X)$ of a
space $X$ leads to two inquiries.

\begin{quote}
\textbf{Problem A.} Specify an additive basis for the graded abelian group
$H^{\ast}(X)$ that encodes the geometric formation of $X$ (e.g. a cell
decomposition of $X$).

\textbf{Problem B.} Determine the table of multiplications between these base elements.
\end{quote}

\bigskip

It is plausible that if $X$ is a flag manifold $G/H$, a uniform solution to
Problem A is afforded by the\textbf{ }\textsl{Basis Theorem}\textbf{ }from the
\textsl{Schubert enumerative calculus} [S$_{2}$]. This was originated by
Ehresmann for the Grassmannians $G_{n,k}$ of $k$-dimensional subspaces in
$\mathbb{C}^{n}$ in 1934 [E], extended to the case where $G$ is a matrix group
by Bruhat in 1954, and completed for all compact connected Lie groups by
Chevalley in 1958 [Ch$_{2}$]. We briefly recall the result.

\bigskip

Let $W$ and $W^{^{\prime}}$ be the Weyl groups of $G$ and $H$ respectively.
The set $W/W^{^{\prime}}$ of left cosets of $W^{^{\prime}}$ in $W$ can be
identified with the subset of $W$:

\begin{center}
$\overline{W}=\{w\in W\mid l(w_{1})\geq l(w)$ for all $w_{1}\in wW^{^{\prime}%
}\}$,\ 
\end{center}

\noindent where $l:W\rightarrow\mathbb{Z}$ is the length function relative to
a fixed maximal torus $T$ in $G$ [BGG, 5.1. Proposition]. The key fact is that
the space $G/H$ admits a canonical decomposition into cells indexed by
elements of $\overline{W}$

\begin{enumerate}
\item[(1.1)] $\qquad\qquad G/H=\underset{w\in\overline{W}}{\cup}X_{w}$,
$\quad\dim X_{w}=2l(w)$,
\end{enumerate}

\noindent with each cell $X_{w}$ the closure of an algebraic affine space,
known as a \textsl{Schubert variety} in $G/H$ [Ch$_{2}$, BGG].

Since only even dimensional cells are involved in the decomposition (1.1), the
set of fundamental classes\textsl{ }$[X_{w}]\in H_{2l(w)}(G/H)$,
$w\in\overline{W}$,\textsl{ }forms an additive basis of $H_{\ast}(G/H)$. The
cocycle class $P_{w}\in H^{2l(w)}(G/H)$, $w\in\overline{W}$, defined by the
Kronecker pairing as $\left\langle P_{w},[X_{u}]\right\rangle =\delta_{w,u}$,
$w,u\in\overline{W}$, is called the \textsl{Schubert class corresponding to
}$w$. The solution to Problem A can be stated in (cf. [BGG])

\begin{quote}
\textbf{Basis Theorem.} \textsl{The set of Schubert classes }$\{P_{w}\mid
$\textsl{\ }$w\in\overline{W}\}$\textsl{\ constitutes an additive basis for
the ring }$H^{\ast}(G/H)$\textsl{.}
\end{quote}

\bigskip

One of the direct consequences of the basis Theorem is that the product of two
arbitrary Schubert classes can be expressed in terms of Schubert classes.
Precisely, given $u,v\in\overline{W}$, one has the expression

\begin{center}
$P_{u}\cdot P_{v}=\sum\limits_{l(w)=l(u)+l(v),w\in\overline{W}}a_{u,v}%
^{w}P_{w}$, $a_{u,v}^{w}\in\mathbb{Z}$
\end{center}

\noindent in $H^{\ast}(G/H)$. Thus, in the case of $X=G/H$, Problem B has a
concrete form.

\begin{quote}
\textbf{Problem B'.} \textsl{Determine the structure constants }$a_{u,v}^{w}$
\textsl{of the ring }$H^{\ast}(G/H)$ \textsl{for }$w,u,v\in\overline{W}%
$\textsl{ with} $l(w)=l(u)+l(v)$\textsl{.}
\end{quote}

\bigskip

Originated in the pioneer works of Schubert on enumerative geometry from 1874
and spurred by Hilbert's fifteenth problem, the study of Problem B' has a long
and outstanding history even for the very special case $G=U(n)$ and
$H=U(k)\times U(n-k)$, where $U(n)$ is the unitary group of rank $n$ (cf.
[K]). The corresponding flag manifold $G/H$ is \textsl{the Grassmannian
}$G_{n,k}$\textsl{ }of $k$-planes through the origin in $\mathbb{C}^{n}$, and
the solution to Problem B' is given by the classical \textsl{Pieri
formula}\footnote{In order to find a formula for the degrees of Schubert
varieties on the Grassmanian, Schubert himself developed a special case of the
Pieri formula [K].} and the \textsl{Littlewood-Richardson rule}%
\footnote{Classically, the Littlewood-Richardson rule describes the
multiplicative rule of Schur symmetric functions. It was first stated by
Littlewood and Richardson in 1934 [LR] and completly proofs appeared in the
70's (see ``Note and references'' in [M, p.148]). Lesieur noticed in 1947 [L]
that the multiplicative rule of Schubert classes in the Grassmanian formly
coincides with that of Schur functions. That is, the Littlewood-Richardson
rule can also be considered as the rule for multiplying Schubert classes in
the Grassmanians.}. We refer to the articles [KL] by Kleiman-Laksov and [St]
by Stanley for full expositions of these results respectively from geometric
approach and from combinatorial view-point.

During the past half century, many achievements have been made in extending
the knowledge on the $a_{u,v}^{w}$ from the $G_{n,k}$ to flag manifolds of
other types. See [Ch$_{1}$], [Mo], [BGG], [D], [LS$_{2}$], [HB], [KK], [Wi],
[BS$_{1}$-BS$_{3}$], [S$_{2}$], [PR$_{1}$- PR$_{3}$], [Bi].

\bigskip

Early in 1953, Borel introduced a method to compute the cohomology algebra
$H^{\ast}(G/H;\mathbb{R})$ (with real coefficients) using spectral sequence
technique [Bo$_{1}$, Bo$_{2}$, B, TW, W]. In the results so obtained the
algebra $H^{\ast}(G/H;\mathbb{R})$ was characterized algebraically in terms of
generators-relations, in which the basis theorem that implies the geometric
structure of the space $G/H$ was absent\footnote{In the intersection theory,
the basis Theorem is important for it guarantees that the rational equivalence
class of a subvariety in $G/H$ can be expressed in term of the base elements
and therefore, the intersection multiplicities of arbitrary subvarities in
$G/H$ can be computed in terms of the $a_{u,v}^{w}$.}. In recent years, in
order to recover from Borel's description of the algebra $H^{\ast
}(G/H;\mathbb{R})$ the polynomial representatives of Schubert classes so that
explicit computation for the $a_{u,v}^{w}$ is possible, various theories
of\textsl{ Schubert polynomials} were developed for the cases where $G$ is a
matrix group and $H\subset G$ is a maximal torus (cf. [S$_{2}$], [LS$_{1}$],
[Be], [BH], [BJS], [FK], [FS], [Fu], [LPR], [Ma]).

\bigskip

Combining the ideas of the Bott-Samelson resolutions of Schubert varieties
[BS, Han] and the enumerative formula on a twisted product of $2$ spheres
developed in [Du$_{2}$], the first author obtained in [Du$_{3}$] a formula
expressing the structure constant $a_{u,v}^{w}$ in terms of Cartan numbers of
$G$. It was also announced in [Du$_{3}$] that, based on the formula, a program
to compute the $a_{u,v}^{w}$ can be compiled. This paper is devoted to explain
the algorithm in details.

Consequently, the algorithm gives a method to compute the integral cohomology
ring $H^{\ast}(G/H)$ independent of the classical spectral sequence method due
to Leray [L$_{1}$,L$_{2}$] and Borel [Bo$_{1}$, Bo$_{2}$]. It has also served
the purpose to indicate our further programs computing Steenrod operations on
$G/H$ [DZ], and multiplying Demazure basis (resp. Grothendieck basis) in the
Grothendieck cohomology of $G/H$ [Du$_{4}$].

\bigskip

The paper is so arranged. In Section 2 we recall the formula for the
$a_{u,v}^{w}$ from [Du$_{3}$]. In Section 3 we resolve Problem B' into two
algorithms entitled ``\textsl{Decompositions}'' and ``\textsl{L-R
coefficients}''. The functions of the algorithms are implemented respectively
in section 4 and 5.

Explicit computation in the cohomology (i.e. the Chow ring) of such classical
spaces as flag varieties is not only required by the effective computability
of problems from enumerative geometry [K], but also related to many problems
from geometry and topology [IM, H, Du$_{1}$]. In order to demonstrate that our
algorithm is effective, samples of computational results from the program are
explained and tabulated in Section 6.

\bigskip

It is worth to mention that there have been excellent codes for multiplying
Schubert classes in the Grassmannians $G_{n,k}$ (cf. the programs SYMMETRICA
at Bayreuth

\begin{center}%
$<$%
http://www. mathe2.uni-bayeuth.de%
$>$%
;
\end{center}

\noindent the programs ACE at Marne La Valle

\begin{center}%
$<$%
http://phalanstere.univ-mlv.fr%
$>$%
,
\end{center}

\noindent and the program LITTLWOOD-RICHARDSON CALCULATOR at Aarhus

\begin{center}%
$<$%
http://home.imf.au.dk/abuch%
$>$%
).
\end{center}

\noindent Instead of being type-specific, our program applies uniformly to all
$G/H$.

\section{The formula}

This section recalls the formula for the $a_{u,v}^{w}$ from [Du$_{3}$]. A few
preliminary notations will be needed. Throughout this paper $G$ is a compact
connected Lie group with a fixed maximal torus $T$. We set $n=\dim T$.

\smallskip\ \ \ 

Equip the Lie algebra $L(G)$ of $G$ with an inner product $(,)$ so that the
adjoint representation acts as isometries of $L(G)$. The\textsl{\ Cartan
subalgebra} of $G$ is the Euclidean subspace $L(T)$ of $L(G)$ [Hu, p.80].

The restriction of the exponential map $\exp:L(G)\rightarrow G$ to $L(T)$
defines a set $D(G)$ of $m=\frac{1}{2}(\dim G-n)$ hyperplanes in $L(T)$, i.e.
the set of\textsl{\ singular hyperplanes }through the origin in $L(T)$. These
planes divide $L(T)$ into finitely many convex cones, called the \textsl{Weyl
chambers} of $G$. The reflections $\sigma$ of $L(T)$ in the these planes
generate \textsl{the Weyl group} $W$ of $G$.

Fix, once and for all, a regular point $\alpha\in L(T)\backslash\underset{L\in
D(G)}{\cup}L$ and let $\Delta=\{\beta_{1},\cdots,\beta_{n}\}$ be the set of
simple roots relative to $\alpha$ [Hu, p.47]. For a $1\leq i\leq n$, write
$\sigma_{i}\in W$ \ for the reflection of $L(T)$ in the singular plane
$L_{\beta_{i}}\in D(G)$ corresponding to the root $\beta_{i}$. The $\sigma
_{i}$ are called \textsl{simple reflections }[Hu, 42].

Recall that for $1\leq i,j\leq n$, \textsl{the Cartan number} $\beta_{i}%
\circ\beta_{j}=:2(\beta_{i},\beta_{j})/(\beta_{j},\beta_{j})$ is always an
integer (only $0,\pm1,\pm2,\pm3$ can occur) [Hu, p.39, p.55].

\bigskip

It is known that the set of simple reflections $\{\sigma_{i}\mid1\leq i\leq
n\}$ generates $W$. That is, any $w\in W$ admits a factorization of the form

\begin{enumerate}
\item[(2.1)] $\qquad\qquad w=\sigma_{i_{1}}\circ\cdots\circ\sigma_{i_{k}}$, .
\end{enumerate}

\begin{quote}
\textbf{Definition 1}. The \textsl{length} $l(w)$ of a $w\in W$ is the least
number of factors in all decompositions of $w$ in the form (2.1). The
decomposition (2.1) is said \textsl{reduced} if $k=l(w)$.

If (2.1) is a reduced decomposition, the $k\times k$ (strictly upper
triangular) matrix $A_{w}=(a_{s,t})$ with

$\qquad\qquad a_{s,t}=\{%
\begin{array}
[c]{c}%
0\text{ if }s\geq t\text{;\qquad}\\
-\beta_{i_{t}}\circ\beta_{i_{s}}\text{ if }s<t
\end{array}
$

is called \textsl{the Cartan matrix of }$w$\textsl{ associated to the
decomposition }(2.1).\ \ \ \ 
\end{quote}

Let $\mathbb{Z}[x_{1},\cdots,x_{k}]=\oplus_{r\geq0}\mathbb{Z}[x_{1}%
,\cdots,x_{k}]^{(r)}$ be the ring of integral polynomials in $x_{1}%
,\cdots,x_{k}$, graded by $\mid x_{i}\mid=1$.

\begin{quote}
\textbf{Definition 2. }Given an $k\times k$ strictly upper triangular integer
matrix $A=(a_{i,j})$, the \textsl{triangular operator} associated to $A$ is
the homomorphism $T_{A}:$ $\mathbb{Z}[x_{1},\cdots,x_{k}]^{(k)}\rightarrow
\mathbb{Z}$ defined recursively by the following \textsl{elimination laws}.

1) if $h\in\mathbb{Z}[x_{1},\cdot\cdot\cdot,x_{k-1}]^{(k)}$, then $T_{A}(h)=0$;

2) if $k=1$ (consequently $A=(0)$), then $T_{A}(x_{1})=1$;

3) if $h\in\mathbb{Z}[x_{1},\cdot\cdot\cdot,x_{k-1}]^{(k-r)}$ with $r\geq1$, then
\end{quote}

\begin{center}
$T_{A}(hx_{k}^{r})=T_{A^{\prime}}(h(a_{1,k}x_{1}+\cdots+a_{k-1,k}%
x_{k-1})^{r-1})$,
\end{center}

\begin{quote}
where $A^{\prime}$ is the ($(k-1)\times(k-1)$ strictly upper triangular)
matrix obtained from $A$ by deleting the $k^{th}$ column and the $k^{th}$ row.
\end{quote}

By additivity, $T_{A}$ is defined for every $f\in\mathbb{Z}[x_{1},\cdots
,x_{k}]^{(k)}$ using the unique expansion $f=\Sigma h_{r}x_{k}^{r}$ with
$h_{r}\in\mathbb{Z}[x_{1},\cdots,x_{k-1}]^{(k-r)}$.

\bigskip

\textbf{Example.} Definition 2 implies an effective algorithm to evaluate
$T_{A}$.

\begin{quote}
For $k=2$ and $A_{1}=\left(
\begin{array}
[c]{cc}%
0 & a\\
0 & 0
\end{array}
\right)  $, then $T_{A_{1}}:$ $\mathbb{Z}[x_{1},x_{2}]^{(2)}\rightarrow
\mathbb{Z}$ is given by
\end{quote}

$\qquad\qquad T_{A_{1}}(x_{1}^{2})=0$,

$\qquad\qquad T_{A_{1}}(x_{1}x_{2})=T_{A_{1}^{\prime}}(x_{1})=1$ and

$\qquad\qquad T_{A_{1}}(x_{2}^{2})=T_{A_{1}^{\prime}}(ax_{1})=a$.

\bigskip

\begin{quote}
For $k=3$ and $A_{2}=\left(
\begin{array}
[c]{ccc}%
0 & a & b\\
0 & 0 & c\\
0 & 0 & 0
\end{array}
\right)  $, then $A_{2}^{\prime}=A_{1}$ and $T_{A_{2}}:$ $\mathbb{Z}%
[x_{1},x_{2},x_{3}]^{(3)}\rightarrow\mathbb{Z}$ is given by
\end{quote}

\begin{center}
$T_{A_{2}}(x_{1}^{r_{1}}x_{2}^{r_{2}}x_{3}^{r_{3}})=\{%
\begin{array}
[c]{c}%
0\text{, if }r_{3}=0\text{ and\qquad\qquad\qquad\qquad\qquad}\\
T_{A_{1}}(x_{1}^{r_{1}}x_{2}^{r_{2}}(bx_{1}+cx_{2})^{r_{3}-1}),\text{ if
}r_{3}\geq1\text{,\ }%
\end{array}
$
\end{center}

\begin{quote}
\noindent where $r_{1}+r_{2}+r_{3}=3$, and where $T_{A_{1}}$ is calculated in
the above.
\end{quote}

\bigskip

Assume that $w=\sigma_{i_{1}}\circ\cdots\circ\sigma_{i_{k}}$, $1\leq
i_{1},\cdots,i_{k}\leq n$, is a reduced decomposition of $w\in\overline{W}$,
and let $A_{w}=(a_{s,t})_{k\times k}$ be the associated Cartan matrix. For a
subset $L=[j_{1},\cdots,j_{r}]\subseteq\lbrack1,\cdots,k]$ we put $\mid
L\mid=r$ and set

\begin{center}
$\sigma_{L}=\sigma_{i_{j_{1}}}\circ\cdots\circ\sigma_{i_{j_{r}}}\in W$;$\quad
x_{L}=x_{j_{1}}\cdots x_{j_{r}}\in\mathbb{Z}[x_{1},\cdots,x_{k}]$.
\end{center}

\noindent The solution to Problem B' is (cf. [Du$_{3}$])

\begin{quote}
\textbf{The formula.} \textsl{If }$u,v\in\overline{W}$\textsl{ with
}$l(w)=l(u)+l(v)$\textsl{, then}

$\qquad\qquad a_{u,v}^{w}=T_{A_{w}}[(\sum\limits_{\mid L\mid=l(u),\text{
}\sigma_{L}=u}x_{L})(\sum\limits_{\mid K\mid=l(v),\text{ }\sigma_{K}=v}%
x_{K})]$\textsl{,}

\textsl{where }$L,K\subseteq\lbrack1,\cdots,k]$\textsl{.}
\end{quote}

\bigskip

The subsequent sections are devoted to clarify the algorithm implicitly
contained in the formula.

\section{The structure of the algorithm}

Let $L(T)$ be the Cartan subalgebra of $G$ and let $\Delta=\{\beta_{1}%
,\cdots,\beta_{n}\}\subset L(T)$ be the set of simple roots of $G$ relative to
the regular point $\alpha\in L(T)$ (cf. Section 2). The \textsl{Cartan matrix}
\ of $G$ is the $n\times n$ integral matrix $C=(c_{ij})_{n\times n}$ defined by

\begin{center}
$c_{ij}=2(\beta_{i},\beta_{j})/(\beta_{j},\beta_{j})$, $1\leq i,j\leq n$.
\end{center}

\noindent It is well known that (cf. [Hu])

\begin{quote}
\textbf{Fact 1.} \textsl{All simply connected compact semi-simple Lie groups
are classified by their Cartan matrices.}
\end{quote}

For a subset $K=[i_{1},\cdots,i_{d}]\subset\lbrack1,\cdots,n]$ let $b\in
L(T)\backslash\{0\}$ be a point lying \textsl{exactly} in the singular
hyperplanes $L_{\beta_{i_{1}}},\cdots,L_{\beta_{i_{d}}}$; namely,

\begin{enumerate}
\item[(3.1)] $\qquad\qquad\qquad b\in\bigcap\limits_{i\in K}$ $L_{\beta_{i}%
}\backslash\bigcup\limits_{j\in J}L_{\beta_{j}}$ ($\in L(T)\backslash
\bigcup\limits_{j\in J}L_{\beta_{j}}$ if $K=\emptyset$)
\end{enumerate}

\noindent where $J$ is the complement of $K$ in $[1,\cdots,n]$. Denote by
$H_{K}$ the centralizer of the $1$-parameter subgroup $\{\exp(tb)\mid
t\in\mathbb{R}\}$ in $G$. It can be shown that (cf. [BHi, 13.5-13.6]))

\begin{quote}
\textbf{Fact 2.} \textsl{The isomorphism type of the Lie group }$H_{K}%
$\textsl{ depends only on the subset }$K$\textsl{ and not on a specific choice
of b in (3.1). Further, every centralizer }$H$ \textsl{of a one-parameter
subgroup in }$G$\textsl{ is conjugate in }$G$\textsl{ to one of the subgroups
}$H_{K}$\textsl{.}
\end{quote}

By Fact 2 we may assume that $H$ is of the form $H_{K}$ for some
$K\subset\lbrack1,\cdots,n]$. Summarizing Fact 1 and 2 we have

\begin{quote}
\textbf{Lemma 1.} \textsl{A complete set of numerical invariants required to
determine a flag manifold }$G/H$\textsl{ consists of}

\textsl{1) a Cartan matrix }$C=(c_{ij})_{n\times n}$\textsl{ (to determine
}$G$\textsl{);}

\textsl{2) a subset }$K=[i_{1},\cdots,i_{d}]\subset\lbrack1,\cdots,n]$
\textsl{(to specify }$H\subset G$\textsl{).}
\end{quote}

The implementation of our program essentially consists of two algorithms,
whose functions may be briefed as follows.

\bigskip

\textbf{Algorithm A.} \textsl{Decompositions.}

\textbf{Input:} \textsl{A Cartan matrix }$C=(c_{ij})_{n\times n}$\textsl{, and
a subset }$K\subset\lbrack1,\ldots,n]$\textsl{.}

\textbf{Output: }\textsl{The coset }$\overline{W}$\textsl{ being presented by
a reduced decomposition for each} $w\in\overline{W}.$

\bigskip

\textbf{Remark 1.} In [Ste, Section 1]\textbf{ }Stembridge described an
algorithm for the problem of \textsl{finding a reduced decomposition for a
given }$w\in W$. This requests less than what Algorithm A concerns.

\bigskip

\textbf{Algorithm B. }\textsl{L-R coefficients}.

\textbf{Input:}\textsl{ }$u$\textsl{, }$v$\textsl{, }$w\in\overline{W}%
$\textsl{ with }$l(u)+l(v)=l(w)$\textsl{.}

\textbf{Output:} $a_{u,v}^{w}\in\mathbb{Z}$.

\bigskip

The details of these algorithms will be given respectively in the coming two sections.

\bigskip

It is clear from the above discussion that our algorithms reduce the structure
constants $a_{u,v}^{w}$ directly to the Cartan matrix $C=(c_{ij})_{n\times n}$
and the subset $K\subset\lbrack1,\ldots,n]$: the simplest and minimum set of
constants \textsl{by which all flag manifolds }$G/H$\textsl{ are classified}
(cf. Lemma 1). Because of this feature it is functional equally for
computations in all $G/H$.

\section{Algorithm A}

We show in \textbf{4.1} the fashion by which the Weyl groups $W^{^{\prime}%
}\subset W$ arise from the Cartan matrix $C=(c_{ij})_{n\times n}$ and the
subset $K\subset\lbrack1,\cdots,n]$. In \textbf{4.2} a numerical
representation for $W$ is introduced. Based on the terminologies developed in
\textbf{4.1} and \textbf{4.2}, Algorithm A is given in \textbf{4.3}.

\bigskip

\textbf{4.1. Constructing the Weyl groups }$W^{^{\prime}}\subset W$
\textbf{from the Cartan matrix. }Let $\Gamma$ be the free $\mathbb{Z}$-module
with $n$ generators $\omega_{1},\cdots,\omega_{n}$, and let $Aut(\Gamma)$ be
the group of automorphisms of $\Gamma$.

Given a Cartan matrix $C=(c_{ij})_{n\times n}$ of a Lie group $G$ with rank
$n$, define $n$ endomorphisms $\sigma_{k}$ of $\Gamma$ (in term of Cartan
numbers) by

\begin{enumerate}
\item[(4.1)] \noindent$\qquad\qquad\sigma_{k}(\omega_{i})=\{%
\begin{array}
[c]{c}%
\omega_{i}\text{ if }k\neq i\text{;\quad\qquad\qquad\quad}\\
\omega_{i}-\Sigma_{1\leq j\leq n}c_{ij}\omega_{j}\text{ if }k=i
\end{array}
$, $1\leq k\leq n$.
\end{enumerate}

\noindent It is straightforward to verify that $\sigma_{k}^{2}=Id$. In
particular, $\sigma_{k}\in Aut(\Gamma)$.

\begin{quote}
\textbf{Lemma 2.} \textsl{The subgroup of }$Aut(\Gamma)$\textsl{\ generated by
}$\sigma_{1},\cdots,\sigma_{n}$\textsl{\ is isomorphic to }$W$\textsl{, the
Weyl group of }$G$\textsl{.}

\textsl{For a subset }$K\subset\lbrack1,\cdots,n]$\textsl{,} \textsl{the
subgroup }$W^{^{\prime}}$ \textsl{of }$W$\textsl{\ generated by }$\{\sigma
_{i}\mid i\in K\}$\textsl{\ is isomorphic the Weyl group of }$H_{K}$ (cf.
Section 3).
\end{quote}

\textbf{Proof }(cf. Proof of Theorem 1 in [DZZ]). Let $t$ be the real vector
space spanned by $\omega_{1},\cdots,\omega_{n}$; namely, $t=\Gamma
\otimes\mathbb{R}$. In term of the Cartan matrix $C=(c_{ij})_{n\times n}$ we
introduce in $t$ the vectors $\beta_{1},\cdots,\beta_{n}$ by

\begin{center}
$\beta_{i}=c_{i1}\omega_{1}+\cdots+c_{in}\omega_{n}$,
\end{center}

\noindent and define an Euclidean metric on $t$ by

\begin{enumerate}
\item[(4.2)] \noindent$\qquad\qquad\qquad2(\beta_{i},\frac{\beta_{j}}%
{(\beta_{j},\beta_{j})})=c_{ij};$ $(\beta_{1},\beta_{1})=1$.
\end{enumerate}

\noindent Then

(a) $t$ can be identified with the Cartan subalgebra $L(T)$ of $G$ under which
the vectors $\beta_{1},\cdots,\beta_{n}$ corresponds to the set $\Delta$ of
simple roots of $G$ (cf. Section 2);

(b) with respect to the metric (4.2), the induced action of $\sigma_{k}$ on
$t=L(T)$ is the reflection in the hyperplane $L_{\beta_{k}}$ perpendicular to
the $\beta_{k}$;

(c) under the identification $t=L(T)$ specified in (a), the basis $\omega
_{1},\cdots,\omega_{n}$ of $\Gamma$ agrees with the set of \textsl{the
fundamental dominant weights} relative to $\Delta$ [Hu, p.67] (Geometrically,
positive multiples of $\omega_{1},\cdots,\omega_{n}$ form the edges of the
Weyl chamber in $L(T)$ corresponding to $\Delta$).

Lemma 2 follows directly from (b) and (c).$\square$

\bigskip

\textbf{4.2.} \textbf{A numerical representation of Weyl groups }In the theory
of Lie algebras the vector $\delta=\omega_{1}+\cdots+\omega_{n}\in
\Gamma\subset$ $t=\Gamma\otimes\mathbb{R}$ is well known as a \textsl{strongly
dominant weight} [Hu, p.30]. For a $w\in W$ consider the expression in $\Gamma$

\begin{center}
$w(\delta)=b_{1}\omega_{1}+\cdots+b_{n}\omega_{n}$, $b_{i}\in\mathbb{Z}$.
\end{center}

\begin{quote}
\textbf{Definition 3.} The correspondence $b:W\rightarrow\mathbb{Z}^{n}$ by
$b(w)=(b_{1},\cdots,b_{n})$ will be called \textsl{the numerical
representation of }$W$.

\textbf{Lemma 3.} \textsl{The numerical representation }$b:W\rightarrow
\mathbb{Z}^{n}$\textsl{ is faithful and satisfies }$b_{i}\neq0$\textsl{ for
all }$w\in W$\textsl{ and }$1\leq i\leq n$\textsl{.}
\end{quote}

\textbf{Proof.} By (c) in the proof of Lemma 2, $\delta\in t$ is a regular
point in the Weyl chamber determined by $\Delta$. Lemma 3 comes from the
geometric fact that the action of the Weyl group $W$ on the orbit of any
regular point is simply transitive.$\square$

\bigskip

The formula (4.1), together with additivity of the $\sigma_{k}$, is sufficient
to compute the coordinates of $b(w)$ from the Cartan numbers and any
decomposition of $w\in W$ into products of the $\sigma_{i}$, as the following
algorithm shows.

\bigskip

\textbf{Algorithm 1.} Computing $b(w)$.

\textbf{Input:} A sequence $1\leq i_{1},\cdots,i_{m}\leq n$.

\textbf{Output:} $b(w)$ for $w=\sigma_{i_{1}}\circ\cdots\circ\sigma_{i_{m}}$.

\textsl{Procedure:} Begin with the sum $p_{0}=\omega_{1}+\cdots+\omega_{n}$.

\textsl{Step 1.} Substituting in $p_{0}$ the term $\omega_{i_{m}}$ by
$\omega_{i_{m}}-\Sigma_{1\leq j\leq n}c_{i_{m}j}\omega_{j}$ to get $p_{1}$;

\textsl{Step 2.} Substituting in $p_{1}$ the term $\omega_{i_{m-1}}$ by
$\omega_{i_{m-1}}-\Sigma_{1\leq j\leq n}c_{i_{m-1}j}\omega_{j}$ to get $p_{2}$;

$\qquad\qquad\vdots$

\textsl{Step m.} Substituting in $p_{m-1}$ the term $\omega_{i_{1}}$ by
$\omega_{i_{1}}-\Sigma_{1\leq j\leq n}c_{i_{1}j}\omega_{j}$ to get $p_{m}$;

\textsl{Step m+1.} If $p_{m}=b_{1}\omega_{1}+\cdots+b_{n}\omega_{n}$ then
$b(w)=(b_{1},\cdots,b_{n})$.

\bigskip

We conclude this subsection with two useful properties of the numerical
representation of a Weyl group given in Definition 3. Let $l:W\rightarrow
\mathbb{Z}$ be the length function on $W$. As in Section 1 we identify
$\overline{W}$ with the subset of $W$

\begin{center}
$\overline{W}=\{w\in W\mid l(w)\leq l(u)$ for all $u\in wW^{^{\prime}}\}$.
\end{center}

\begin{quote}
\textbf{Lemma 4.} \textsl{Let }$w\in W$\textsl{ be with }$b(w)=(b_{1}%
,\cdots,b_{n})$\textsl{ and }$b(w^{-1})=(\overline{b}_{1},\cdots,\overline
{b}_{n})$\textsl{. Then}

\textsl{(i) }$l(\sigma_{i}w)=l(w)-1$\textsl{ if and only if }$b_{i}<0$\textsl{;}

\textsl{(ii) }$w\in\overline{W}$\textsl{ if and only if }$\overline{b}_{i}%
>0$\textsl{ for all }$i\in K$\textsl{.}
\end{quote}

\textbf{Proof.} The metric on $L(T)$ yields the relations

\begin{enumerate}
\item[(4.3)] $\qquad\qquad(\omega_{i},\beta_{j}/(\beta_{j},\beta_{j}%
))=\delta_{ij}$
\end{enumerate}

\noindent between the simple roots $\beta_{j}$ and the corresponding
fundamental dominant weights $\omega_{i}$ [Hu, p.67]. By [BGG, 2.3 Corollary],
$l(\sigma_{i}w)=l(w)-1$ if and only if $(w(\delta),\beta_{i})<0$. The latter
is equivalent to $b_{i}<0$ in views of (4.3) and $w(\delta)=b_{1}\omega
_{1}+\cdots+b_{n}\omega_{n}$. This verifies (i).

Similarly, assertion (ii) follows from the following alternative description
for $\overline{W}$ (cf. [BGG, 5.1. Proposition, (iii)])

\begin{center}
$\overline{W}=\{w\in W\mid(w^{-1}(\delta),\beta_{i})>0$\textsl{ for all }$i\in
K\}$.$\square$
\end{center}

\bigskip

\textbf{4.3. Construction of the coset} $\overline{W}=W/W^{\prime}$. Let
$l:W\rightarrow\mathbb{Z}$ be the length function on $W$. We put $\overline
{W}^{k}=\{w\in\overline{W}\mid l(w)=k\}$, $k=0,1,2,\cdots$. Then, as is clear,
$\overline{W}=\coprod\limits_{k\geq0}\overline{W}^{k}$. The problem concerned
by Algorithm A may be reduced to

\begin{quote}
\textbf{Problem C.} \textsl{Enumerate elements in }$\overline{W}^{k}$\textsl{
(i.e. in }$\overline{W}$\textsl{), }$k\geq0$\textsl{, by their reduced decompositions.}
\end{quote}

Before presenting Algorithm A (i.e. the solution to Problem C) we note that

\begin{enumerate}
\item[(4.4)] If the set $\overline{W}^{k}$ is given in term of certain reduced
decompositions of its elements, then $\overline{W}^{k}$ becomes an ordered set
with the order specified by

$\qquad\qquad\sigma_{i_{1}}\circ\cdots\circ\sigma_{i_{k}}<\sigma_{j_{1}}%
\circ\cdots\circ\sigma_{j_{k}}$

if there exists some $s\leq k$ such that $i_{t}=j_{t}$ for all $t<s$ \ but
$i_{s}<j_{s}$.

\item[(4.5)] If $X$ and $Y$ are two ordered sets, then the product $X\times Y$
is furnished with the canonical order as:

$``(x,y)<(x^{\prime},y^{\prime})$ if and only if $x<x^{\prime}$ or
$x=x^{\prime}$ but $y<y^{\prime}"$.
\end{enumerate}

\bigskip

The solution for Problem C is known when $k=0,1$

\begin{center}
$\overline{W}^{0}=\{id\}$; $\quad\overline{W}^{1}=\{\sigma_{j}\mid j\in J\}$,
\end{center}

\noindent where $id$ is the identity of $W$ and where $J$ is the complement of
$K$ in $[1,\cdots,n]$. In general, Algorithm A enables one to build up
$\overline{W}^{k}$ from $\overline{W}^{k-1}$.

\bigskip

\textbf{Algorithm A. }\textsl{Decompositions.}

\textbf{Input.} The set $\overline{W}^{k-1}$ being presented by certain
reduced decompositions of its elements.

\textbf{Output.} The set $\overline{W}^{k}$ being presented by certain reduced
decompositions of its elements.

\textsl{Procedure:} Set $V=\{1,\cdots,n\}\times\overline{W}^{k-1}$. Repeat the
following steps for all elements in $V$ in accordance with the order on $V$
(cf. (4.5)). Begin with empty sets $S=\emptyset$, $R=\emptyset$.

\textsl{Step 1.} For a $v=(i,\sigma_{i_{1}}\circ\cdots\circ\sigma_{i_{k-1}%
})\in V$ form the product $w=\sigma_{i}\circ\sigma_{i_{1}}\circ\cdots
\circ\sigma_{i_{k-1}}$.

\textsl{Step 2.} Call Algorithm 1 to obtain

\begin{center}
$b(w)=(b_{1},\cdots,b_{n})$ and $b(w^{-1})=(\overline{b}_{1},\cdots
,\overline{b}_{n})$;
\end{center}

\textsl{Step 3.} If \quad1) $b_{i}<0;$

\qquad\qquad\qquad2) $\overline{b}_{i}>0$ for all $i\in K;$

\qquad\qquad\qquad3) $(b_{1},\cdots,b_{n})\notin R$,

\noindent add $\sigma_{i}\circ\sigma_{i_{1}}\circ\cdots\circ\sigma_{i_{k-1}}$
to $S$; add $b(w)=(b_{1},\cdots,b_{n})$ to $R$;

The program terminates at $S=\overline{W}^{k}$.

\bigskip

\textbf{Explanation}. We verify the last clause in Algorithm A. Firstly, Lemma
7 in [DZZ] claims that any $w\in\overline{W}^{k}$ admits a decomposition
$w=\sigma_{i}\circ\sigma_{i_{1}}\circ\cdots\circ\sigma_{i_{k-1}}$ for some
$(i,\sigma_{i_{1}}\circ\cdots\circ\sigma_{i_{k-1}})\in V$. This explains the
role the set $V$ plays in the algorithm. Next, the first two conditions in
Step 3 guarantees that $\sigma_{i}\circ\sigma_{i_{1}}\circ\cdots\circ
\sigma_{i_{k-1}}\in\overline{W}^{k}$ by Lemma 4. Finally, the third constraint
in Step 3 rejects a second reduced decomposition of some $w\in\overline{W}%
^{k}$ being included in $\overline{W}^{k}$ (by Lemma 3).$\square$

\bigskip

\textbf{Remark 2.} If $K=\emptyset$, then $\overline{W}=W$ (the whole group).
In this case $H=T$ (a maximal torus in $G$) and Step 2 and 3 in Algorithm A
can be simplified as

\textsl{Step 2.} Call Algorithm 1 to obtain $b(w)=(b_{1},\cdots,b_{n})$;

\textsl{Step 3.} If $b_{i}<0$ and if $(b_{1},\cdots,b_{n})\notin R$, add
$\sigma_{i}\circ\sigma_{i_{1}}\circ\cdots\circ\sigma_{i_{k-1}}$ to $S$; add
$b(w)=(b_{1},\cdots,b_{n})$ to $R$.

\bigskip

\textbf{Remark 3.} Based on the word representation of Weyl groups, a
different program solving Problem C was given in [DZZ]. In comparison, the use
of the numerical representation simplifies the presentation of Algorithm A.

\section{Algorithm B}

Algorithm A presents us the coset $\overline{W}=\coprod\limits_{k\geq
0}\overline{W}^{k}$ by certain reduced decomposition of its elements. Based on
this we explain \textsl{L-R coefficients},\textsl{ }the algorithm computing
$a_{u,v}^{w}$.

By the notion $L\subset\lbrack1,\cdots,k],\left|  L\right|  =r$, we mean that
$L$ is a sequence $(j_{1},\cdots,j_{r})$ of $r$ integers satisfying

\begin{center}
$1\leq j_{1}<\cdots<j_{r}\leq k$.
\end{center}

\noindent For two integers $1\leq r\leq k$ let the set $V(k,r)=\{L\mid
L\subset\lbrack1,\cdots,k],\left|  L\right|  =r\}$ be equipped with the
obvious ordering (cf. (4.3)).

For a $w=\sigma_{i_{1}}\circ\cdots\circ\sigma_{i_{k}}\in\overline{W}$ and a
$u\in\overline{W}$ with $l(u)=r<k$, we set

\begin{enumerate}
\item[(5.1)] \noindent$\qquad\qquad\qquad p_{w}(u)=\sum\limits_{L\in
V(k,r),\sigma_{L}=u}x_{L}\in\mathbb{Z}[x_{1},\cdots,x_{k}]^{(r)}$,
\end{enumerate}

\noindent where $\sigma_{L}=\sigma_{i_{j_{1}}}\circ\cdots\circ\sigma
_{i_{j_{r}}}$ if $L=[j_{1},\cdots,j_{r}]$. Using these notations our formula
(cf. Section 2) can be simplified as

\begin{enumerate}
\item[(5.2)] $\qquad\qquad\qquad a_{u,v}^{w}=T_{A_{w}}[p_{w}(u)p_{w}(v)]$.
\end{enumerate}

\bigskip

We begin by pointing out that (5.1) suggests the following algorithm
specifying the polynomial $p_{w}(u)$.

\bigskip

\textbf{Algorithm 2.} Computing $p_{w}(u)\in\mathbb{Z}[x_{1},\cdots
,x_{k}]^{(r)}$.

\textbf{Input:} $w=\sigma_{i_{1}}\circ\cdots\circ\sigma_{i_{k}}\in\overline
{W}^{k}$and $u\in\overline{W}^{r}$ with $b(u)=(b_{1},\cdots,b_{n})$.

\textbf{Output:} $p_{w}(u)$.

\textsl{Procedure:} Repeat the following steps for all $L\in V(k,r)$ in
accordance with the order on $V(k,r)$. Initiate the polynomial $p=p(x_{1}%
,\cdots,x_{k})$ as zero.

\textsl{Step 1.} For a $L\in V(k,r)$ call algorithm 1 to get $b(\sigma_{L})$;

\textsl{Step 2.} If $b(\sigma_{L})=b(u)$ add $x_{L}$ to $p$.

The program terminates at $p=p_{w}(u)$.

\bigskip

If $A=(a_{ij})_{k\times k}$ is matrix of rank $k$ and if $1\leq r\leq k-1$,
then the notion $(a_{ij})_{r\times r}$ clearly stands for the matrix of rank
$r$ obtained from $A$ by deleting the last $(k-r)$ rows and columns.

Let $A=(a_{ij})_{k\times k}$ be a strictly upper triangular integral matrix of
rank $k$. Consider the triangular operator $T_{A}:$ $\mathbb{Z}[x_{1}%
,\cdots,x_{k}]^{(k)}\rightarrow\mathbb{Z}$ given in Definition 2.

\bigskip

\textbf{Algorithm 3.} Computing $T_{A}:$ $\mathbb{Z}[x_{1},\cdots,x_{k}%
]^{(k)}\rightarrow\mathbb{Z}$.

\textbf{Input:} A strictly upper triangular integral matrix $A=(a_{ij}%
)_{k\times k}$ and $\ $a polynomial $p=p(x_{1},\cdots,x_{k})\in\mathbb{Z}%
[x_{1},\cdots,x_{k}]^{(k)}$

\textbf{Output:} $T_{A}(p)\in\mathbb{Z}$.

\textsl{Procedure:}\textbf{ }Recursion.

\textsl{Step 1.} Express $p$ as a polynomial in $x_{k}$; i.e.

\begin{center}
$p=h_{0}+h_{1}x_{k}+\sum\limits_{2\leq r\leq k}h_{r}x_{k}^{r}$, $h_{r}%
\in\mathbb{Z}[x_{1},\cdots,x_{k-1}]^{(k-r)}$,
\end{center}

\noindent and set

\begin{center}
$p_{1}=h_{1}+\sum\limits_{2\leq r\leq k}h_{r}(a_{1,k}x_{1}+\cdots
+a_{k-1,k}x_{k-1})^{r-1}$($\in\mathbb{Z}[x_{1},\cdots,x_{k-1}]^{(k-1)}$).
\end{center}

\textsl{Step 2.} Repeat step 1 for $A_{1}=(a_{ij})_{(k-1)\times(k-1)}$ and
$p=p_{1}$ to get $p_{2}\in\mathbb{Z}[x_{1},\cdots,x_{k-2}]^{(k-2)}$.

$\qquad\qquad\vdots$

\textsl{Step k+1.} If $p_{k}=a\in\mathbb{Z}$, then $T_{A}(p)=a$.

\bigskip

\textbf{Algorithm B. }\textsl{L-R coefficients.}

\textbf{Input:} $w=\sigma_{i_{1}}\circ\cdots\circ\sigma_{i_{k}}\in\overline
{W}^{k}$, $(u,v)\in\overline{W}^{r}\times\overline{W}^{k-r}$

\textbf{Output: }$a_{u,v}^{w}\in\mathbb{Z}$.

\textsl{Procedure:} Let $A_{w}$ \ be the Cartan matrix of $w$ related to the
decomposition (it can be read directly from the Cartan matrix of $G$ and the
decomposition $w=\sigma_{i_{1}}\circ\cdots\circ\sigma_{i_{k}}$. cf. Definition 1).

\textsl{Step 1.} Call algorithm 2 to get $p_{w}(u)$ and $p_{w}(v)$;

\textsl{Step 2.} Call algorithm 3 to get $T_{A_{w}}(p_{w}(u)\cdot p_{w}(v))$.

\textsl{Step 3.} If $T_{A_{w}}(p_{w}(u)\cdot p_{w}(v))=a$, then $a_{u,v}%
^{w}=a$ (by (5.2)).

\bigskip

\textbf{Remark 4.} Based on Algorithm B, a parallel program to expand the product

\begin{center}
$P_{u}\cdot P_{v}=\sum\limits_{w\in\overline{W}^{k}}a_{u,v}^{w}P_{w}$
\end{center}

\noindent for given $(u,v)\in\overline{W}^{r}\times\overline{W}^{k-r}$ can be
easily implemented. The order on $\overline{W}^{k}$ can be employed to assign
each $w\in\overline{W}^{k}$ a computing unit.

\section{Computational examples}

Our algorithm is ready to apply to computation in flag manifolds. Recall that
all compact connected semi-simple irreducible Lie groups fall into four
infinite sequences of matrix groups

\begin{center}
$SU(n)$;$\quad SO(2n)$;$\quad SO(2n+1)$;$\quad Sp(n)$,
\end{center}

\noindent as well as the five exceptional ones

\begin{center}
$G_{2}$,\noindent\qquad\noindent$F_{4}$,$\qquad E_{6}$,$\qquad E_{7}$,$\qquad
E_{8}$.
\end{center}

\noindent The flag manifolds associated to matrix groups have been studied
extensively during the past decades. Here, we choose to work with certain flag
manifolds related to the exceptional Lie groups $E_{n}$, $n=6,7,8$.

\bigskip

Fix a maximal torus $T^{n}\subset E_{n}$ and let $W(n)$ be the Weyl group of
$E_{n}$. Then

\begin{enumerate}
\item[(6.1)] $\qquad\left|  W(n)\right|  =\{%
\begin{array}
[c]{c}%
2^{7}3^{4}5\ \text{\textsl{if} }n=6\text{;}\\
2^{10}3^{4}57\ \text{\textsl{if} }n=7\text{;}\\
2^{14}3^{5}5^{2}7\ \text{\textsl{if} }n=8\text{;}%
\end{array}
\quad\dim_{\mathbb{R}}E_{n}=\{%
\begin{array}
[c]{c}%
78\ \text{\textsl{if} }n=6\text{\textsl{;}}\\
133\ \text{\textsl{if} }n=7\text{\textsl{;}}\\
248\ \text{\textsl{if} }n=8\text{\textsl{,}}%
\end{array}
$
\end{enumerate}

\noindent where $\left|  A\right|  $ stands for the cardinality of the set
$A$. Assume that the set of simple roots $\Delta=\{\beta_{1},\cdots,\beta
_{n}\}$ of $E_{n}$ is given and ordered as the vertices of the Dynkin diagram
of $E_{n}$ pictured in [Hu, p.58], and let $K\subset\{1,2,\cdots,n\}$ be the
subset whose complement is $\{2\}$. We have the following information on the
subgroup group $H_{K}$.

\begin{quote}
(a) \textsl{the semisimple part of the subgroup }$H_{K}\subset E_{n}$\textsl{
is }$SU(n)$\textsl{, the special unitary group of order }$n$\textsl{;}

(b) $H_{K}$\textsl{ admits a factorization into the semi-product }$H_{K}%
=S^{1}\cdot SU(n)$\textsl{, where }$S^{1}$\textsl{ is a circle subgroup of the
maximal torus }$T^{n}$\textsl{ in }$E_{n}$\textsl{;}

(c) \textsl{if }$W^{\prime}(n)\subset W(n)$\textsl{ is the Weyl group of
}$H_{K}$\textsl{, then }$\left|  W^{\prime}(n)\right|  =n!$.
\end{quote}

Consequently, if one write $\overline{W}(n)$ for the coset $W^{^{\prime}}(n)$
in $W(n)$, one has

\begin{enumerate}
\item[(6.2)] $\qquad\left|  \overline{W}(n)\right|  =\{%
\begin{array}
[c]{c}%
2^{3}3^{2}\ \text{\textsl{if} }n=6\text{;}\\
2^{6}3^{2}\ \text{\textsl{if} }n=7\text{;}\\
2^{7}3^{3}5\ \text{\textsl{if} }n=8\text{.}%
\end{array}
$\qquad$\dim_{\mathbb{R}}E_{n}/H_{K}=\{%
\begin{array}
[c]{c}%
42\ \text{\textsl{if} }n=6\text{\textsl{;}}\\
84\ \text{\textsl{if} }n=7\text{\textsl{;}}\\
194\ \text{\textsl{if} }n=8\text{\textsl{.}}%
\end{array}
$
\end{enumerate}

\noindent Geometrically, $\overline{W}(n)$ parameterizes Schubert classes on
$E_{n}/H_{K}$ (i.e. the Basis Theorem).

\bigskip

The subset of $\overline{W}(n)$ consisting of elements with length $r$ is
denoted $\overline{W}^{r}(n)$. By (4.4), if the $\overline{W}^{r}(n)$ is
presented by its elements each with a reduced decompositions, then it
naturally becomes an ordered set and therefore, can be alternatively presented as

\begin{enumerate}
\item[(6.3)] $\qquad\qquad\overline{W}^{r}(n)=\{w_{r,i}\mid1\leq i\leq\left|
\overline{W}^{r}\right|  \}$.
\end{enumerate}

\bigskip

In table A$_{n}$ below, we present elements of $\overline{W}(n)$ with length
$r\leq10$ both in terms of their reduced decompositions produced by Algorithm
A, and the index system (6.3) imposed by the decompositions.

The index (6.3) on $\overline{W}^{r}(n)$ is useful in simplifying the
presentation of the intersection multiplicities $a_{u,v}^{w}$. By resorting to
this index system we list in table B$_{n}$ ( $n=6,7,8$) all the $a_{u,v}^{w}$
with $l(w)=9$ and $10$ produced by Algorithm B.

\bigskip

\begin{center}
\textbf{Table A}$_{6}$\textbf{.}{ \textsl{Reduced decomposition of elements in
}$\overline{W}(6)$\textsl{ with length}$\leq10$}\newline
\begin{tabular}
[c]{|l|l|l|}\hline
$w_{i,j}$ & decomposition & \\\hline
$w_{1,1}$ & $\sigma_{2}$ & \\
$w_{2,1}$ & $\sigma_{4}\sigma_{2}$ & \\
$w_{3,1}$ & $\sigma_{3}\sigma_{4}\sigma_{2}$ & \\
$w_{3,2}$ & $\sigma_{5}\sigma_{4}\sigma_{2}$ & \\
$w_{4,1}$ & $\sigma_{1}\sigma_{3}\sigma_{4}\sigma_{2}$ & \\
$w_{4,2}$ & $\sigma_{3}\sigma_{5}\sigma_{4}\sigma_{2}$ & \\
$w_{4,3}$ & $\sigma_{6}\sigma_{5}\sigma_{4}\sigma_{2}$ & \\
$w_{5,1}$ & $\sigma_{1}\sigma_{3}\sigma_{5}\sigma_{4}\sigma_{2}$ & \\
$w_{5,2}$ & $\sigma_{3}\sigma_{6}\sigma_{5}\sigma_{4}\sigma_{2}$ & \\
$w_{5,3}$ & $\sigma_{4}\sigma_{3}\sigma_{5}\sigma_{4}\sigma_{2}$ & \\
$w_{6,1}$ & $\sigma_{1}\sigma_{3}\sigma_{6}\sigma_{5}\sigma_{4}\sigma_{2}$ &
\\
$w_{6,2}$ & $\sigma_{1}\sigma_{4}\sigma_{3}\sigma_{5}\sigma_{4}\sigma_{2}$ &
\\
$w_{6,3}$ & $\sigma_{2}\sigma_{4}\sigma_{3}\sigma_{5}\sigma_{4}\sigma_{2}$ &
\\
$w_{6,4}$ & $\sigma_{4}\sigma_{3}\sigma_{6}\sigma_{5}\sigma_{4}\sigma_{2}$ &
\\
$w_{7,1}$ & $\sigma_{1}\sigma_{2}\sigma_{4}\sigma_{3}\sigma_{5}\sigma
_{4}\sigma_{2}$ & \\
$w_{7,2}$ & $\sigma_{1}\sigma_{4}\sigma_{3}\sigma_{6}\sigma_{5}\sigma
_{4}\sigma_{2}$ & \\
$w_{7,3}$ & $\sigma_{2}\sigma_{4}\sigma_{3}\sigma_{6}\sigma_{5}\sigma
_{4}\sigma_{2}$ & \\
$w_{7,4}$ & $\sigma_{3}\sigma_{1}\sigma_{4}\sigma_{3}\sigma_{5}\sigma
_{4}\sigma_{2}$ & \\\hline
\end{tabular}
\negthinspace%
\begin{tabular}
[c]{|l|l|l|}\hline
$w_{i,j}$ & decomposition & \\\hline
$w_{7,5}$ & $\sigma_{5}\sigma_{4}\sigma_{3}\sigma_{6}\sigma_{5}\sigma
_{4}\sigma_{2}$ & \\
$w_{8,1}$ & $\sigma_{1}\sigma_{2}\sigma_{4}\sigma_{3}\sigma_{6}\sigma
_{5}\sigma_{4}\sigma_{2}$ & \\
$w_{8,2}$ & $\sigma_{1}\sigma_{5}\sigma_{4}\sigma_{3}\sigma_{6}\sigma
_{5}\sigma_{4}\sigma_{2}$ & \\
$w_{8,3}$ & $\sigma_{2}\sigma_{3}\sigma_{1}\sigma_{4}\sigma_{3}\sigma
_{5}\sigma_{4}\sigma_{2}$ & \\
$w_{8,4}$ & $\sigma_{2}\sigma_{5}\sigma_{4}\sigma_{3}\sigma_{6}\sigma
_{5}\sigma_{4}\sigma_{2}$ & \\
$w_{8,5}$ & $\sigma_{3}\sigma_{1}\sigma_{4}\sigma_{3}\sigma_{6}\sigma
_{5}\sigma_{4}\sigma_{2}$ & \\
$w_{9,1}$ & $\sigma_{1}\sigma_{2}\sigma_{5}\sigma_{4}\sigma_{3}\sigma
_{6}\sigma_{5}\sigma_{4}\sigma_{2}$ & \\
$w_{9,2}$ & $\sigma_{2}\sigma_{3}\sigma_{1}\sigma_{4}\sigma_{3}\sigma
_{6}\sigma_{5}\sigma_{4}\sigma_{2}$ & \\
$w_{9,3}$ & $\sigma_{3}\sigma_{1}\sigma_{5}\sigma_{4}\sigma_{3}\sigma
_{6}\sigma_{5}\sigma_{4}\sigma_{2}$ & \\
$w_{9,4}$ & $\sigma_{4}\sigma_{2}\sigma_{3}\sigma_{1}\sigma_{4}\sigma
_{3}\sigma_{5}\sigma_{4}\sigma_{2}$ & \\
$w_{9,5}$ & $\sigma_{4}\sigma_{2}\sigma_{5}\sigma_{4}\sigma_{3}\sigma
_{6}\sigma_{5}\sigma_{4}\sigma_{2}$ & \\
$w_{10,1}$ & $\sigma_{1}\sigma_{4}\sigma_{2}\sigma_{5}\sigma_{4}\sigma
_{3}\sigma_{6}\sigma_{5}\sigma_{4}\sigma_{2}$ & \\
$w_{10,2}$ & $\sigma_{2}\sigma_{3}\sigma_{1}\sigma_{5}\sigma_{4}\sigma
_{3}\sigma_{6}\sigma_{5}\sigma_{4}\sigma_{2}$ & \\
$w_{10,3}$ & $\sigma_{3}\sigma_{4}\sigma_{2}\sigma_{5}\sigma_{4}\sigma
_{3}\sigma_{6}\sigma_{5}\sigma_{4}\sigma_{2}$ & \\
$w_{10,4}$ & $\sigma_{4}\sigma_{2}\sigma_{3}\sigma_{1}\sigma_{4}\sigma
_{3}\sigma_{6}\sigma_{5}\sigma_{4}\sigma_{2}$ & \\
$w_{10,5}$ & $\sigma_{4}\sigma_{3}\sigma_{1}\sigma_{5}\sigma_{4}\sigma
_{3}\sigma_{6}\sigma_{5}\sigma_{4}\sigma_{2}$ & \\
$w_{10,6}$ & $\sigma_{5}\sigma_{4}\sigma_{2}\sigma_{3}\sigma_{1}\sigma
_{4}\sigma_{3}\sigma_{5}\sigma_{4}\sigma_{2}$ & \\
&  & \\\hline
\end{tabular}
\end{center}

\newpage\renewcommand{\arraystretch}{0.95}

\begin{center}
{\setlength{\tabcolsep}{1mm} }
\end{center}

\textbf{Table A}$_{7}${\small . \textsl{Reduced decomposition of elements in
}$\overline{W}(7)$\textsl{ with length}$\leq10$\textsl{\newline }%
\begin{tabular}
[c]{|l|l|}\hline
$w_{i,j}$ & decomposition\\\hline
$w_{1,1}$ & $\sigma_{2}$\\
$w_{2,1}$ & $\sigma_{4}\sigma_{2}$\\
$w_{3,1}$ & $\sigma_{3}\sigma_{4}\sigma_{2}$\\
$w_{3,2}$ & $\sigma_{5}\sigma_{4}\sigma_{2}$\\
$w_{4,1}$ & $\sigma_{1}\sigma_{3}\sigma_{4}\sigma_{2}$\\
$w_{4,2}$ & $\sigma_{3}\sigma_{5}\sigma_{4}\sigma_{2}$\\
$w_{4,3}$ & $\sigma_{6}\sigma_{5}\sigma_{4}\sigma_{2}$\\
$w_{5,1}$ & $\sigma_{1}\sigma_{3}\sigma_{5}\sigma_{4}\sigma_{2}$\\
$w_{5,2}$ & $\sigma_{3}\sigma_{6}\sigma_{5}\sigma_{4}\sigma_{2}$\\
$w_{5,3}$ & $\sigma_{4}\sigma_{3}\sigma_{5}\sigma_{4}\sigma_{2}$\\
$w_{5,4}$ & $\sigma_{7}\sigma_{6}\sigma_{5}\sigma_{4}\sigma_{2}$\\
$w_{6,1}$ & $\sigma_{1}\sigma_{3}\sigma_{6}\sigma_{5}\sigma_{4}\sigma_{2}$\\
$w_{6,2}$ & $\sigma_{1}\sigma_{4}\sigma_{3}\sigma_{5}\sigma_{4}\sigma_{2}$\\
$w_{6,3}$ & $\sigma_{2}\sigma_{4}\sigma_{3}\sigma_{5}\sigma_{4}\sigma_{2}$\\
$w_{6,4}$ & $\sigma_{3}\sigma_{7}\sigma_{6}\sigma_{5}\sigma_{4}\sigma_{2}$\\
$w_{6,5}$ & $\sigma_{4}\sigma_{3}\sigma_{6}\sigma_{5}\sigma_{4}\sigma_{2}$\\
$w_{7,1}$ & $\sigma_{1}\sigma_{2}\sigma_{4}\sigma_{3}\sigma_{5}\sigma
_{4}\sigma_{2}$\\
$w_{7,2}$ & $\sigma_{1}\sigma_{3}\sigma_{7}\sigma_{6}\sigma_{5}\sigma
_{4}\sigma_{2}$\\\hline
\end{tabular}
\negthinspace\negthinspace\negthinspace\
\begin{tabular}
[c]{|l|l|l|}\hline
$w_{i,j}$ & decomposition & \\\hline
$w_{7,3}$ & $\sigma_{1}\sigma_{4}\sigma_{3}\sigma_{6}\sigma_{5}\sigma
_{4}\sigma_{2}$ & \\
$w_{7,4}$ & $\sigma_{2}\sigma_{4}\sigma_{3}\sigma_{6}\sigma_{5}\sigma
_{4}\sigma_{2}$ & \\
$w_{7,5}$ & $\sigma_{3}\sigma_{1}\sigma_{4}\sigma_{3}\sigma_{5}\sigma
_{4}\sigma_{2}$ & \\
$w_{7,6}$ & $\sigma_{4}\sigma_{3}\sigma_{7}\sigma_{6}\sigma_{5}\sigma
_{4}\sigma_{2}$ & \\
$w_{7,7}$ & $\sigma_{5}\sigma_{4}\sigma_{3}\sigma_{6}\sigma_{5}\sigma
_{4}\sigma_{2}$ & \\
$w_{8,1}$ & $\sigma_{1}\sigma_{2}\sigma_{4}\sigma_{3}\sigma_{6}\sigma
_{5}\sigma_{4}\sigma_{2}$ & \\
$w_{8,2}$ & $\sigma_{1}\sigma_{4}\sigma_{3}\sigma_{7}\sigma_{6}\sigma
_{5}\sigma_{4}\sigma_{2}$ & \\
$w_{8,3}$ & $\sigma_{1}\sigma_{5}\sigma_{4}\sigma_{3}\sigma_{6}\sigma
_{5}\sigma_{4}\sigma_{2}$ & \\
$w_{8,4}$ & $\sigma_{2}\sigma_{3}\sigma_{1}\sigma_{4}\sigma_{3}\sigma
_{5}\sigma_{4}\sigma_{2}$ & \\
$w_{8,5}$ & $\sigma_{2}\sigma_{4}\sigma_{3}\sigma_{7}\sigma_{6}\sigma
_{5}\sigma_{4}\sigma_{2}$ & \\
$w_{8,6}$ & $\sigma_{2}\sigma_{5}\sigma_{4}\sigma_{3}\sigma_{6}\sigma
_{5}\sigma_{4}\sigma_{2}$ & \\
$w_{8,7}$ & $\sigma_{3}\sigma_{1}\sigma_{4}\sigma_{3}\sigma_{6}\sigma
_{5}\sigma_{4}\sigma_{2}$ & \\
$w_{8,8}$ & $\sigma_{5}\sigma_{4}\sigma_{3}\sigma_{7}\sigma_{6}\sigma
_{5}\sigma_{4}\sigma_{2}$ & \\
$w_{9,1}$ & $\sigma_{1}\sigma_{2}\sigma_{4}\sigma_{3}\sigma_{7}\sigma
_{6}\sigma_{5}\sigma_{4}\sigma_{2}$ & \\
$w_{9,2}$ & $\sigma_{1}\sigma_{2}\sigma_{5}\sigma_{4}\sigma_{3}\sigma
_{6}\sigma_{5}\sigma_{4}\sigma_{2}$ & \\
$w_{9,3}$ & $\sigma_{1}\sigma_{5}\sigma_{4}\sigma_{3}\sigma_{7}\sigma
_{6}\sigma_{5}\sigma_{4}\sigma_{2}$ & \\
$w_{9,4}$ & $\sigma_{2}\sigma_{3}\sigma_{1}\sigma_{4}\sigma_{3}\sigma
_{6}\sigma_{5}\sigma_{4}\sigma_{2}$ & \\
$w_{9,5}$ & $\sigma_{2}\sigma_{5}\sigma_{4}\sigma_{3}\sigma_{7}\sigma
_{6}\sigma_{5}\sigma_{4}\sigma_{2}$ & \\\hline
\end{tabular}
\negthinspace\negthinspace\negthinspace\
\begin{tabular}
[c]{|l|l|l|}\hline
$w_{i,j}$ & decomposition & \\\hline
$w_{9,6}$ & $\sigma_{3}\sigma_{1}\sigma_{4}\sigma_{3}\sigma_{7}\sigma
_{6}\sigma_{5}\sigma_{4}\sigma_{2}$ & \\
$w_{9,7}$ & $\sigma_{3}\sigma_{1}\sigma_{5}\sigma_{4}\sigma_{3}\sigma
_{6}\sigma_{5}\sigma_{4}\sigma_{2}$ & \\
$w_{9,8}$ & $\sigma_{4}\sigma_{2}\sigma_{3}\sigma_{1}\sigma_{4}\sigma
_{3}\sigma_{5}\sigma_{4}\sigma_{2}$ & \\
$w_{9,9}$ & $\sigma_{4}\sigma_{2}\sigma_{5}\sigma_{4}\sigma_{3}\sigma
_{6}\sigma_{5}\sigma_{4}\sigma_{2}$ & \\
$w_{9,10}$ & $\sigma_{6}\sigma_{5}\sigma_{4}\sigma_{3}\sigma_{7}\sigma
_{6}\sigma_{5}\sigma_{4}\sigma_{2}$ & \\
$w_{10,1}$ & $\sigma_{1}\sigma_{2}\sigma_{5}\sigma_{4}\sigma_{3}\sigma
_{7}\sigma_{6}\sigma_{5}\sigma_{4}\sigma_{2}$ & \\
$w_{10,2}$ & $\sigma_{1}\sigma_{4}\sigma_{2}\sigma_{5}\sigma_{4}\sigma
_{3}\sigma_{6}\sigma_{5}\sigma_{4}\sigma_{2}$ & \\
$w_{10,3}$ & $\sigma_{1}\sigma_{6}\sigma_{5}\sigma_{4}\sigma_{3}\sigma
_{7}\sigma_{6}\sigma_{5}\sigma_{4}\sigma_{2}$ & \\
$w_{10,4}$ & $\sigma_{2}\sigma_{3}\sigma_{1}\sigma_{4}\sigma_{3}\sigma
_{7}\sigma_{6}\sigma_{5}\sigma_{4}\sigma_{2}$ & \\
$w_{10,5}$ & $\sigma_{2}\sigma_{3}\sigma_{1}\sigma_{5}\sigma_{4}\sigma
_{3}\sigma_{6}\sigma_{5}\sigma_{4}\sigma_{2}$ & \\
$w_{10,6}$ & $\sigma_{2}\sigma_{6}\sigma_{5}\sigma_{4}\sigma_{3}\sigma
_{7}\sigma_{6}\sigma_{5}\sigma_{4}\sigma_{2}$ & \\
$w_{10,7}$ & $\sigma_{3}\sigma_{1}\sigma_{5}\sigma_{4}\sigma_{3}\sigma
_{7}\sigma_{6}\sigma_{5}\sigma_{4}\sigma_{2}$ & \\
$w_{10,8}$ & $\sigma_{3}\sigma_{4}\sigma_{2}\sigma_{5}\sigma_{4}\sigma
_{3}\sigma_{6}\sigma_{5}\sigma_{4}\sigma_{2}$ & \\
$w_{10,9}$ & $\sigma_{4}\sigma_{2}\sigma_{3}\sigma_{1}\sigma_{4}\sigma
_{3}\sigma_{6}\sigma_{5}\sigma_{4}\sigma_{2}$ & \\
$w_{10,10}$ & $\sigma_{4}\sigma_{2}\sigma_{5}\sigma_{4}\sigma_{3}\sigma
_{7}\sigma_{6}\sigma_{5}\sigma_{4}\sigma_{2}$ & \\
$w_{10,11}$ & $\sigma_{4}\sigma_{3}\sigma_{1}\sigma_{5}\sigma_{4}\sigma
_{3}\sigma_{6}\sigma_{5}\sigma_{4}\sigma_{2}$ & \\
$w_{10,12}$ & $\sigma_{5}\sigma_{4}\sigma_{2}\sigma_{3}\sigma_{1}\sigma
_{4}\sigma_{3}\sigma_{5}\sigma_{4}\sigma_{2}$ & \\
&  & \\\hline
\end{tabular}
}

\begin{center}
{ \setlength{\tabcolsep}{0.5mm} \textbf{Table A}$_{8}$.{\small \textsl{Reduced
decomposition of elements in }$\overline{W}(8)$\textsl{ with length}$\leq
10$.\newline
\begin{tabular}
[c]{|l|l|l|}\hline
$w_{i,j}$ & decomposition & \\\hline
$w_{1,1}$ & $\sigma_{2}$ & \\
$w_{2,1}$ & $\sigma_{4}\sigma_{2}$ & \\
$w_{3,1}$ & $\sigma_{3}\sigma_{4}\sigma_{2}$ & \\
$w_{3,2}$ & $\sigma_{5}\sigma_{4}\sigma_{2}$ & \\
$w_{4,1}$ & $\sigma_{1}\sigma_{3}\sigma_{4}\sigma_{2}$ & \\
$w_{4,2}$ & $\sigma_{3}\sigma_{5}\sigma_{4}\sigma_{2}$ & \\
$w_{4,3}$ & $\sigma_{6}\sigma_{5}\sigma_{4}\sigma_{2}$ & \\
$w_{5,1}$ & $\sigma_{1}\sigma_{3}\sigma_{5}\sigma_{4}\sigma_{2}$ & \\
$w_{5,2}$ & $\sigma_{3}\sigma_{6}\sigma_{5}\sigma_{4}\sigma_{2}$ & \\
$w_{5,3}$ & $\sigma_{4}\sigma_{3}\sigma_{5}\sigma_{4}\sigma_{2}$ & \\
$w_{5,4}$ & $\sigma_{7}\sigma_{6}\sigma_{5}\sigma_{4}\sigma_{2}$ & \\
$w_{6,1}$ & $\sigma_{1}\sigma_{3}\sigma_{6}\sigma_{5}\sigma_{4}\sigma_{2}$ &
\\
$w_{6,2}$ & $\sigma_{1}\sigma_{4}\sigma_{3}\sigma_{5}\sigma_{4}\sigma_{2}$ &
\\
$w_{6,3}$ & $\sigma_{2}\sigma_{4}\sigma_{3}\sigma_{5}\sigma_{4}\sigma_{2}$ &
\\
$w_{6,4}$ & $\sigma_{3}\sigma_{7}\sigma_{6}\sigma_{5}\sigma_{4}\sigma_{2}$ &
\\
$w_{6,5}$ & $\sigma_{4}\sigma_{3}\sigma_{6}\sigma_{5}\sigma_{4}\sigma_{2}$ &
\\
$w_{6,6}$ & $\sigma_{8}\sigma_{7}\sigma_{6}\sigma_{5}\sigma_{4}\sigma_{2}$ &
\\
$w_{7,1}$ & $\sigma_{1}\sigma_{2}\sigma_{4}\sigma_{3}\sigma_{5}\sigma
_{4}\sigma_{2}$ & \\
$w_{7,2}$ & $\sigma_{1}\sigma_{3}\sigma_{7}\sigma_{6}\sigma_{5}\sigma
_{4}\sigma_{2}$ & \\
$w_{7,3}$ & $\sigma_{1}\sigma_{4}\sigma_{3}\sigma_{6}\sigma_{5}\sigma
_{4}\sigma_{2}$ & \\
$w_{7,4}$ & $\sigma_{2}\sigma_{4}\sigma_{3}\sigma_{6}\sigma_{5}\sigma
_{4}\sigma_{2}$ & \\
$w_{7,5}$ & $\sigma_{3}\sigma_{1}\sigma_{4}\sigma_{3}\sigma_{5}\sigma
_{4}\sigma_{2}$ & \\\hline
\end{tabular}
\negthinspace\negthinspace\negthinspace\
\begin{tabular}
[c]{|l|l|l|}\hline
$w_{i,j}$ & decomposition & \\\hline
$w_{7,6}$ & $\sigma_{3}\sigma_{8}\sigma_{7}\sigma_{6}\sigma_{5}\sigma
_{4}\sigma_{2}$ & \\
$w_{7,7}$ & $\sigma_{4}\sigma_{3}\sigma_{7}\sigma_{6}\sigma_{5}\sigma
_{4}\sigma_{2}$ & \\
$w_{7,8}$ & $\sigma_{5}\sigma_{4}\sigma_{3}\sigma_{6}\sigma_{5}\sigma
_{4}\sigma_{2}$ & \\
$w_{8,1}$ & $\sigma_{1}\sigma_{2}\sigma_{4}\sigma_{3}\sigma_{6}\sigma
_{5}\sigma_{4}\sigma_{2}$ & \\
$w_{8,2}$ & $\sigma_{1}\sigma_{3}\sigma_{8}\sigma_{7}\sigma_{6}\sigma
_{5}\sigma_{4}\sigma_{2}$ & \\
$w_{8,3}$ & $\sigma_{1}\sigma_{4}\sigma_{3}\sigma_{7}\sigma_{6}\sigma
_{5}\sigma_{4}\sigma_{2}$ & \\
$w_{8,4}$ & $\sigma_{1}\sigma_{5}\sigma_{4}\sigma_{3}\sigma_{6}\sigma
_{5}\sigma_{4}\sigma_{2}$ & \\
$w_{8,5}$ & $\sigma_{2}\sigma_{3}\sigma_{1}\sigma_{4}\sigma_{3}\sigma
_{5}\sigma_{4}\sigma_{2}$ & \\
$w_{8,6}$ & $\sigma_{2}\sigma_{4}\sigma_{3}\sigma_{7}\sigma_{6}\sigma
_{5}\sigma_{4}\sigma_{2}$ & \\
$w_{8,7}$ & $\sigma_{2}\sigma_{5}\sigma_{4}\sigma_{3}\sigma_{6}\sigma
_{5}\sigma_{4}\sigma_{2}$ & \\
$w_{8,8}$ & $\sigma_{3}\sigma_{1}\sigma_{4}\sigma_{3}\sigma_{6}\sigma
_{5}\sigma_{4}\sigma_{2}$ & \\
$w_{8,9}$ & $\sigma_{4}\sigma_{3}\sigma_{8}\sigma_{7}\sigma_{6}\sigma
_{5}\sigma_{4}\sigma_{2}$ & \\
$w_{8,10}$ & $\sigma_{5}\sigma_{4}\sigma_{3}\sigma_{7}\sigma_{6}\sigma
_{5}\sigma_{4}\sigma_{2}$ & \\
$w_{9,1}$ & $\sigma_{1}\sigma_{2}\sigma_{4}\sigma_{3}\sigma_{7}\sigma
_{6}\sigma_{5}\sigma_{4}\sigma_{2}$ & \\
$w_{9,2}$ & $\sigma_{1}\sigma_{2}\sigma_{5}\sigma_{4}\sigma_{3}\sigma
_{6}\sigma_{5}\sigma_{4}\sigma_{2}$ & \\
$w_{9,3}$ & $\sigma_{1}\sigma_{4}\sigma_{3}\sigma_{8}\sigma_{7}\sigma
_{6}\sigma_{5}\sigma_{4}\sigma_{2}$ & \\
$w_{9,4}$ & $\sigma_{1}\sigma_{5}\sigma_{4}\sigma_{3}\sigma_{7}\sigma
_{6}\sigma_{5}\sigma_{4}\sigma_{2}$ & \\
$w_{9,5}$ & $\sigma_{2}\sigma_{3}\sigma_{1}\sigma_{4}\sigma_{3}\sigma
_{6}\sigma_{5}\sigma_{4}\sigma_{2}$ & \\
$w_{9,6}$ & $\sigma_{2}\sigma_{4}\sigma_{3}\sigma_{8}\sigma_{7}\sigma
_{6}\sigma_{5}\sigma_{4}\sigma_{2}$ & \\
$w_{9,7}$ & $\sigma_{2}\sigma_{5}\sigma_{4}\sigma_{3}\sigma_{7}\sigma
_{6}\sigma_{5}\sigma_{4}\sigma_{2}$ & \\
$w_{9,8}$ & $\sigma_{3}\sigma_{1}\sigma_{4}\sigma_{3}\sigma_{7}\sigma
_{6}\sigma_{5}\sigma_{4}\sigma_{2}$ & \\
$w_{9,9}$ & $\sigma_{3}\sigma_{1}\sigma_{5}\sigma_{4}\sigma_{3}\sigma
_{6}\sigma_{5}\sigma_{4}\sigma_{2}$ & \\\hline
\end{tabular}
\negthinspace\negthinspace\negthinspace\
\begin{tabular}
[c]{|l|l|l|}\hline
$w_{i,j}$ & decomposition & \\\hline
$w_{9,10}$ & $\sigma_{4}\sigma_{2}\sigma_{3}\sigma_{1}\sigma_{4}\sigma
_{3}\sigma_{5}\sigma_{4}\sigma_{2}$ & \\
$w_{9,11}$ & $\sigma_{4}\sigma_{2}\sigma_{5}\sigma_{4}\sigma_{3}\sigma
_{6}\sigma_{5}\sigma_{4}\sigma_{2}$ & \\
$w_{9,12}$ & $\sigma_{5}\sigma_{4}\sigma_{3}\sigma_{8}\sigma_{7}\sigma
_{6}\sigma_{5}\sigma_{4}\sigma_{2}$ & \\
$w_{9,13}$ & $\sigma_{6}\sigma_{5}\sigma_{4}\sigma_{3}\sigma_{7}\sigma
_{6}\sigma_{5}\sigma_{4}\sigma_{2}$ & \\
$w_{10,1}$ & $\sigma_{1}\sigma_{2}\sigma_{4}\sigma_{3}\sigma_{8}\sigma
_{7}\sigma_{6}\sigma_{5}\sigma_{4}\sigma_{2}$ & \\
$w_{10,2}$ & $\sigma_{1}\sigma_{2}\sigma_{5}\sigma_{4}\sigma_{3}\sigma
_{7}\sigma_{6}\sigma_{5}\sigma_{4}\sigma_{2}$ & \\
$w_{10,3}$ & $\sigma_{1}\sigma_{4}\sigma_{2}\sigma_{5}\sigma_{4}\sigma
_{3}\sigma_{6}\sigma_{5}\sigma_{4}\sigma_{2}$ & \\
$w_{10,4}$ & $\sigma_{1}\sigma_{5}\sigma_{4}\sigma_{3}\sigma_{8}\sigma
_{7}\sigma_{6}\sigma_{5}\sigma_{4}\sigma_{2}$ & \\
$w_{10,5}$ & $\sigma_{1}\sigma_{6}\sigma_{5}\sigma_{4}\sigma_{3}\sigma
_{7}\sigma_{6}\sigma_{5}\sigma_{4}\sigma_{2}$ & \\
$w_{10,6}$ & $\sigma_{2}\sigma_{3}\sigma_{1}\sigma_{4}\sigma_{3}\sigma
_{7}\sigma_{6}\sigma_{5}\sigma_{4}\sigma_{2}$ & \\
$w_{10,7}$ & $\sigma_{2}\sigma_{3}\sigma_{1}\sigma_{5}\sigma_{4}\sigma
_{3}\sigma_{6}\sigma_{5}\sigma_{4}\sigma_{2}$ & \\
$w_{10,8}$ & $\sigma_{2}\sigma_{5}\sigma_{4}\sigma_{3}\sigma_{8}\sigma
_{7}\sigma_{6}\sigma_{5}\sigma_{4}\sigma_{2}$ & \\
$w_{10,9}$ & $\sigma_{2}\sigma_{6}\sigma_{5}\sigma_{4}\sigma_{3}\sigma
_{7}\sigma_{6}\sigma_{5}\sigma_{4}\sigma_{2}$ & \\
$w_{10,10}$ & $\sigma_{3}\sigma_{1}\sigma_{4}\sigma_{3}\sigma_{8}\sigma
_{7}\sigma_{6}\sigma_{5}\sigma_{4}\sigma_{2}$ & \\
$w_{10,11}$ & $\sigma_{3}\sigma_{1}\sigma_{5}\sigma_{4}\sigma_{3}\sigma
_{7}\sigma_{6}\sigma_{5}\sigma_{4}\sigma_{2}$ & \\
$w_{10,12}$ & $\sigma_{3}\sigma_{4}\sigma_{2}\sigma_{5}\sigma_{4}\sigma
_{3}\sigma_{6}\sigma_{5}\sigma_{4}\sigma_{2}$ & \\
$w_{10,13}$ & $\sigma_{4}\sigma_{2}\sigma_{3}\sigma_{1}\sigma_{4}\sigma
_{3}\sigma_{6}\sigma_{5}\sigma_{4}\sigma_{2}$ & \\
$w_{10,14}$ & $\sigma_{4}\sigma_{2}\sigma_{5}\sigma_{4}\sigma_{3}\sigma
_{7}\sigma_{6}\sigma_{5}\sigma_{4}\sigma_{2}$ & \\
$w_{10,15}$ & $\sigma_{4}\sigma_{3}\sigma_{1}\sigma_{5}\sigma_{4}\sigma
_{3}\sigma_{6}\sigma_{5}\sigma_{4}\sigma_{2}$ & \\
$w_{10,16}$ & $\sigma_{5}\sigma_{4}\sigma_{2}\sigma_{3}\sigma_{1}\sigma
_{4}\sigma_{3}\sigma_{5}\sigma_{4}\sigma_{2}$ & \\
$w_{10,17}$ & $\sigma_{6}\sigma_{5}\sigma_{4}\sigma_{3}\sigma_{8}\sigma
_{7}\sigma_{6}\sigma_{5}\sigma_{4}\sigma_{2}$ & \\
&  & \\\hline
\end{tabular}
}}
\end{center}

\newpage

{ \setlength{\tabcolsep}{1 mm} \textbf{Table B}$_{6}$. }\textsl{L-R
coefficients for }$E_{6}/S^{1}\cdot SU(6)$\textsl{ }

\begin{center}
{\footnotesize
\begin{tabular}
[c]{|r|r|r|r|r|r|r|}\hline
$u$ & $v$ & \multicolumn{5}{c|}{$w\in W^{9}(6)$}\\\cline{3-7}
&  & $w_{9,1}$ & $w_{9,2}$ & $w_{9,3}$ & $w_{9,4}$ & \ $w_{9,5}$\\\hline
$w_{1,1}$ & $w_{8,1}$ & $1$ & $1$ & $0$ & $0$ & $0$\\
$w_{1,1}$ & $w_{8,2}$ & $1$ & $0$ & $1$ & $0$ & $0$\\
$w_{1,1}$ & $w_{8,3}$ & $0$ & $1$ & $0$ & $1$ & $0$\\
$w_{1,1}$ & $w_{8,4}$ & $1$ & $0$ & $0$ & $0$ & $1$\\
$w_{1,1}$ & $w_{8,5}$ & $0$ & $1$ & $1$ & $0$ & $0$\\
$w_{2,1}$ & $w_{7,1}$ & $1$ & $2$ & $0$ & $1$ & $0$\\
$w_{2,1}$ & $w_{7,2}$ & $2$ & $2$ & $2$ & $0$ & $0$\\
$w_{2,1}$ & $w_{7,3}$ & $2$ & $1$ & $0$ & $0$ & $1$\\
$w_{2,1}$ & $w_{7,4}$ & $0$ & $2$ & $1$ & $1$ & $0$\\
$w_{2,1}$ & $w_{7,5}$ & $2$ & $0$ & $1$ & $0$ & $1$\\
$w_{3,1}$ & $w_{6,1}$ & $1$ & $1$ & $1$ & $0$ & $0$\\
$w_{3,1}$ & $w_{6,2}$ & $1$ & $3$ & $2$ & $1$ & $0$\\
$w_{3,1}$ & $w_{6,3}$ & $2$ & $1$ & $0$ & $1$ & $0$\\
$w_{3,1}$ & $w_{6,4}$ & $3$ & $2$ & $1$ & $0$ & $1$\\
$w_{3,2}$ & $w_{6,1}$ & $1$ & $1$ & $1$ & $0$ & $0$\\
$w_{3,2}$ & $w_{6,2}$ & $2$ & $3$ & $1$ & $1$ & $0$\\
$w_{3,2}$ & $w_{6,3}$ & $1$ & $2$ & $0$ & $0$ & $1$\\
$w_{3,2}$ & $w_{6,4}$ & $3$ & $1$ & $2$ & $0$ & $1$\\
$w_{4,1}$ & $w_{5,1}$ & $0$ & $1$ & $1$ & $0$ & $0$\\
$w_{4,1}$ & $w_{5,2}$ & $1$ & $1$ & $0$ & $0$ & $0$\\
$w_{4,1}$ & $w_{5,3}$ & $1$ & $1$ & $1$ & $1$ & $0$\\
$w_{4,2}$ & $w_{5,1}$ & $2$ & $3$ & $2$ & $1$ & $0$\\
$w_{4,2}$ & $w_{5,2}$ & $3$ & $2$ & $2$ & $0$ & $1$\\
$w_{4,2}$ & $w_{5,3}$ & $5$ & $5$ & $2$ & $1$ & $1$\\
$w_{4,3}$ & $w_{5,1}$ & $1$ & $1$ & $0$ & $0$ & $0$\\
$w_{4,3}$ & $w_{5,2}$ & $1$ & $0$ & $1$ & $0$ & $0$\\
$w_{4,3}$ & $w_{5,3}$ & $1$ & $1$ & $1$ & $0$ & $1$\\\hline
\end{tabular}
}

{\footnotesize
\begin{tabular}
[c]{|r|r|r|r|r|r|r|r|}\hline
$u$ & $v$ & \multicolumn{6}{c|}{$w\in W^{10}(6)$}\\\cline{3-8}
&  & $w_{10,1}$ & $w_{10,2}$ & $w_{10,3}$ & \ $w_{10,4}$ & $w_{10,5}$ &
$w_{10,6}$\\\hline
$w_{1,1}$ & $w_{9,1}$ & $1$ & $1$ & $0$ & $0$ & $0$ & $0$\\
$w_{1,1}$ & $w_{9,2}$ & $0$ & $1$ & $0$ & $1$ & $0$ & $0$\\
$w_{1,1}$ & $w_{9,3}$ & $0$ & $1$ & $0$ & $0$ & $1$ & $0$\\
$w_{1,1}$ & $w_{9,4}$ & $0$ & $0$ & $0$ & $1$ & $0$ & $1$\\
$w_{1,1}$ & $w_{9,5}$ & $1$ & $0$ & $1$ & $0$ & $0$ & $0$\\
$w_{2,1}$ & $w_{8,1}$ & $1$ & $2$ & $0$ & $1$ & $0$ & $0$\\
$w_{2,1}$ & $w_{8,2}$ & $1$ & $2$ & $0$ & $0$ & $1$ & $0$\\
$w_{2,1}$ & $w_{8,3}$ & $0$ & $1$ & $0$ & $2$ & $0$ & $1$\\
$w_{2,1}$ & $w_{8,4}$ & $2$ & $1$ & $1$ & $0$ & $0$ & $0$\\
$w_{2,1}$ & $w_{8,5}$ & $0$ & $2$ & $0$ & $1$ & $1$ & $0$\\
$w_{3,1}$ & $w_{7,1}$ & $0$ & $2$ & $0$ & $1$ & $0$ & $1$\\
$w_{3,1}$ & $w_{7,2}$ & $1$ & $3$ & $0$ & $1$ & $1$ & $0$\\
$w_{3,1}$ & $w_{7,3}$ & $2$ & $1$ & $0$ & $1$ & $0$ & $0$\\
$w_{3,1}$ & $w_{7,4}$ & $0$ & $1$ & $0$ & $2$ & $1$ & $0$\\
$w_{3,1}$ & $w_{7,5}$ & $1$ & $2$ & $1$ & $0$ & $0$ & $0$\\
$w_{3,2}$ & $w_{7,1}$ & $1$ & $1$ & $0$ & $2$ & $0$ & $0$\\
$w_{3,2}$ & $w_{7,2}$ & $1$ & $3$ & $0$ & $1$ & $1$ & $0$\\
$w_{3,2}$ & $w_{7,3}$ & $1$ & $2$ & $1$ & $0$ & $0$ & $0$\\
$w_{3,2}$ & $w_{7,4}$ & $0$ & $2$ & $0$ & $1$ & $0$ & $1$\\
$w_{3,2}$ & $w_{7,5}$ & $2$ & $1$ & $0$ & $0$ & $1$ & $0$\\
$w_{4,1}$ & $w_{6,1}$ & $0$ & $1$ & $0$ & $0$ & $0$ & $0$\\
$w_{4,1}$ & $w_{6,2}$ & $0$ & $1$ & $0$ & $1$ & $1$ & $0$\\
$w_{4,1}$ & $w_{6,3}$ & $0$ & $1$ & $0$ & $0$ & $0$ & $1$\\
$w_{4,1}$ & $w_{6,4}$ & $1$ & $1$ & $0$ & $1$ & $0$ & $0$\\
$w_{4,2}$ & $w_{6,1}$ & $1$ & $2$ & $0$ & $1$ & $1$ & $0$\\
$w_{4,2}$ & $w_{6,2}$ & $1$ & $5$ & $0$ & $3$ & $1$ & $1$\\
$w_{4,2}$ & $w_{6,3}$ & $2$ & $2$ & $0$ & $2$ & $0$ & $0$\\
$w_{4,2}$ & $w_{6,4}$ & $3$ & $5$ & $1$ & $1$ & $1$ & $0$\\
$w_{4,3}$ & $w_{6,1}$ & $0$ & $1$ & $0$ & $0$ & $0$ & $0$\\
$w_{4,3}$ & $w_{6,2}$ & $1$ & $1$ & $0$ & $1$ & $0$ & $0$\\
$w_{4,3}$ & $w_{6,3}$ & $0$ & $1$ & $1$ & $0$ & $0$ & $0$\\
$w_{4,3}$ & $w_{6,4}$ & $1$ & $1$ & $0$ & $0$ & $1$ & $0$\\
$w_{5,1}$ & $w_{5,1}$ & $0$ & $2$ & $0$ & $1$ & $1$ & $0$\\
$w_{5,1}$ & $w_{5,2}$ & $1$ & $2$ & $0$ & $1$ & $0$ & $0$\\
$w_{5,1}$ & $w_{5,3}$ & $1$ & $3$ & $0$ & $2$ & $1$ & $1$\\
$w_{5,2}$ & $w_{5,2}$ & $1$ & $2$ & $0$ & $0$ & $1$ & $0$\\
$w_{5,2}$ & $w_{5,3}$ & $2$ & $3$ & $1$ & $1$ & $1$ & $0$\\
$w_{5,3}$ & $w_{5,3}$ & $3$ & $6$ & $0$ & $3$ & $0$ & $0$\\\hline
\end{tabular}
}
\end{center}

\newpage

\begin{center}
{ \setlength{\tabcolsep}{1.5 mm}
}
\end{center}

{\textbf{Table B}}$_{7}$. \textsl{L-R coefficients for }$E_{7}/S^{1}\cdot SU(7)$

{%
\begin{tabular}
[c]{|r|r|r|r|r|r|r|r|r|r|r|r|}\hline
$u$ & $v$ & \multicolumn{10}{c|}{$w\in W^{9}(7)$}\\\cline{3-12}
&  & $w_{9,1}$ & $w_{9,2}$ & $w_{9,3}$ & $w_{9,4}$ & \ $w_{9,5}$ & $w_{9,6}$ &
$w_{9,7}$ & $w_{9,8}$ & \ $w_{9,9}$ & $w_{9,10}$\\\hline
$w_{1,1}$ & $w_{8,1}$ & $1$ & $1$ & $0$ & $1$ & $0$ & $0$ & $0$ & $0$ & $0$ &
$0$\\
$w_{1,1}$ & $w_{8,2}$ & $1$ & $0$ & $1$ & $0$ & $0$ & $1$ & $0$ & $0$ & $0$ &
$0$\\
$w_{1,1}$ & $w_{8,3}$ & $0$ & $1$ & $1$ & $0$ & $0$ & $0$ & $1$ & $0$ & $0$ &
$0$\\
$w_{1,1}$ & $w_{8,4}$ & $0$ & $0$ & $0$ & $1$ & $0$ & $0$ & $0$ & $1$ & $0$ &
$0$\\
$w_{1,1}$ & $w_{8,5}$ & $1$ & $0$ & $0$ & $0$ & $1$ & $0$ & $0$ & $0$ & $0$ &
$0$\\
$w_{1,1}$ & $w_{8,6}$ & $0$ & $1$ & $0$ & $0$ & $1$ & $0$ & $0$ & $0$ & $1$ &
$0$\\
$w_{1,1}$ & $w_{8,7}$ & $0$ & $0$ & $0$ & $1$ & $0$ & $1$ & $1$ & $0$ & $0$ &
$0$\\
$w_{1,1}$ & $w_{8,8}$ & $0$ & $0$ & $1$ & $0$ & $1$ & $0$ & $0$ & $0$ & $0$ &
$1$\\
$w_{2,1}$ & $w_{7,1}$ & $1$ & $1$ & $0$ & $2$ & $0$ & $0$ & $0$ & $1$ & $0$ &
$0$\\
$w_{2,1}$ & $w_{7,2}$ & $1$ & $0$ & $1$ & $0$ & $0$ & $1$ & $0$ & $0$ & $0$ &
$0$\\
$w_{2,1}$ & $w_{7,3}$ & $2$ & $2$ & $2$ & $2$ & $0$ & $2$ & $2$ & $0$ & $0$ &
$0$\\
$w_{2,1}$ & $w_{7,4}$ & $2$ & $2$ & $0$ & $1$ & $2$ & $0$ & $0$ & $0$ & $1$ &
$0$\\
$w_{2,1}$ & $w_{7,5}$ & $0$ & $0$ & $0$ & $2$ & $0$ & $1$ & $1$ & $1$ & $0$ &
$0$\\
$w_{2,1}$ & $w_{7,6}$ & $2$ & $0$ & $2$ & $0$ & $2$ & $1$ & $0$ & $0$ & $0$ &
$1$\\
$w_{2,1}$ & $w_{7,7}$ & $0$ & $2$ & $2$ & $0$ & $2$ & $0$ & $1$ & $0$ & $1$ &
$1$\\
$w_{3,1}$ & $w_{6,1}$ & $1$ & $1$ & $1$ & $1$ & $0$ & $2$ & $1$ & $0$ & $0$ &
$0$\\
$w_{3,1}$ & $w_{6,2}$ & $1$ & $1$ & $1$ & $3$ & $0$ & $1$ & $2$ & $1$ & $0$ &
$0$\\
$w_{3,1}$ & $w_{6,3}$ & $1$ & $2$ & $0$ & $1$ & $1$ & $0$ & $0$ & $1$ & $0$ &
$0$\\
$w_{3,1}$ & $w_{6,4}$ & $2$ & $0$ & $1$ & $0$ & $1$ & $1$ & $0$ & $0$ & $0$ &
$0$\\
$w_{3,1}$ & $w_{6,5}$ & $3$ & $3$ & $3$ & $2$ & $2$ & $1$ & $1$ & $0$ & $1$ &
$1$\\
$w_{3,2}$ & $w_{6,1}$ & $2$ & $1$ & $2$ & $1$ & $0$ & $1$ & $1$ & $0$ & $0$ &
$0$\\
$w_{3,2}$ & $w_{6,2}$ & $2$ & $2$ & $1$ & $3$ & $0$ & $2$ & $1$ & $1$ & $0$ &
$0$\\
$w_{3,2}$ & $w_{6,3}$ & $2$ & $1$ & $0$ & $2$ & $1$ & $0$ & $0$ & $0$ & $1$ &
$0$\\
$w_{3,2}$ & $w_{6,4}$ & $1$ & $0$ & $2$ & $0$ & $1$ & $1$ & $0$ & $0$ & $0$ &
$1$\\
$w_{3,2}$ & $w_{6,5}$ & $3$ & $3$ & $3$ & $1$ & $4$ & $2$ & $2$ & $0$ & $1$ &
$1$\\
$w_{4,1}$ & $w_{5,1}$ & $0$ & $0$ & $0$ & $1$ & $0$ & $1$ & $1$ & $0$ & $0$ &
$0$\\
$w_{4,1}$ & $w_{5,2}$ & $1$ & $1$ & $1$ & $1$ & $0$ & $1$ & $0$ & $0$ & $0$ &
$0$\\
$w_{4,1}$ & $w_{5,3}$ & $1$ & $1$ & $1$ & $1$ & $0$ & $0$ & $1$ & $1$ & $0$ &
$0$\\
$w_{4,1}$ & $w_{5,4}$ & $1$ & $0$ & $0$ & $0$ & $0$ & $0$ & $0$ & $0$ & $0$ &
$0$\\
$w_{4,2}$ & $w_{5,1}$ & $2$ & $2$ & $2$ & $3$ & $0$ & $2$ & $2$ & $1$ & $0$ &
$0$\\
$w_{4,2}$ & $w_{5,2}$ & $5$ & $3$ & $4$ & $2$ & $3$ & $3$ & $2$ & $0$ & $1$ &
$1$\\
$w_{4,2}$ & $w_{5,3}$ & $4$ & $5$ & $3$ & $5$ & $3$ & $2$ & $2$ & $1$ & $1$ &
$1$\\
$w_{4,2}$ & $w_{5,4}$ & $1$ & $0$ & $1$ & $0$ & $1$ & $1$ & $0$ & $0$ & $0$ &
$0$\\
$w_{4,3}$ & $w_{5,1}$ & $2$ & $1$ & $1$ & $1$ & $0$ & $1$ & $0$ & $0$ & $0$ &
$0$\\
$w_{4,3}$ & $w_{5,2}$ & $1$ & $1$ & $3$ & $0$ & $2$ & $1$ & $1$ & $0$ & $0$ &
$1$\\
$w_{4,3}$ & $w_{5,3}$ & $3$ & $1$ & $1$ & $1$ & $2$ & $2$ & $1$ & $0$ & $1$ &
$0$\\
$w_{4,3}$ & $w_{5,4}$ & $0$ & $0$ & $1$ & $0$ & $0$ & $0$ & $0$ & $0$ & $0$ &
$1$\\\hline
\end{tabular}
}

\begin{center}
\renewcommand{\arraystretch}{0.8} { \setlength{\tabcolsep}{0 mm}
{\footnotesize
}}
\end{center}

{\footnotesize \renewcommand{\arraystretch}{0.8} \setlength{\tabcolsep}{0.5
mm}
\begin{tabular}
[c]{|r|r|r|r|r|r|r|r|r|r|r|r|r|r|}\hline
$u$ & $v$ & \multicolumn{12}{c|}{$w\in W^{10}(7)$}\\\cline{3-14}
&  & $w_{10,1}$ & $w_{10,2}$ & $w_{10,3}$ & \ $w_{10,4}$ & $w_{10,5}$ &
$w_{10,6}$ & $w_{10,7}$ & \ $w_{10,8}$ & $w_{10,9}$ & $w_{10,10}$ &
$w_{10,11}$ & $w_{10,12}$\\\hline
{\small $w_{1,1}$} & $w_{9,1}$ & $1$ & $0$ & $0$ & $1$ & $0$ & $0$ & $0$ & $0$
& $0$ & $0$ & $0$ \  & $0$\\
$w_{1,1}$ & $w_{9,2}$ & $1$ & $1$ & $0$ & $0$ & $1$ & $0$ & $0$ & $0$ & $0$ &
$0$ & $0$ \  & $0$\\
$w_{1,1}$ & $w_{9,3}$ & $1$ & $0$ & $1$ & $0$ & $0$ & $0$ & $1$ & $0$ & $0$ &
$0$ & $0$ \  & $0$\\
$w_{1,1}$ & $w_{9,4}$ & $0$ & $0$ & $0$ & $1$ & $1$ & $0$ & $0$ & $0$ & $1$ &
$0$ & $0$ \  & $0$\\
$w_{1,1}$ & $w_{9,5}$ & $1$ & $0$ & $0$ & $0$ & $0$ & $1$ & $0$ & $0$ & $0$ &
$1$ & $0$ \  & $0$\\
$w_{1,1}$ & $w_{9,6}$ & $0$ & $0$ & $0$ & $1$ & $0$ & $0$ & $1$ & $0$ & $0$ &
$0$ & $0$ \  & $0$\\
$w_{1,1}$ & $w_{9,7}$ & $0$ & $0$ & $0$ & $0$ & $1$ & $0$ & $1$ & $0$ & $0$ &
$0$ & $1$ \  & $0$\\
$w_{1,1}$ & $w_{9,8}$ & $0$ & $0$ & $0$ & $0$ & $0$ & $0$ & $0$ & $0$ & $1$ &
$0$ & $0$ \  & $1$\\
$w_{1,1}$ & $w_{9,9}$ & $0$ & $1$ & $0$ & $0$ & $0$ & $0$ & $0$ & $1$ & $0$ &
$1$ & $0$ \  & $0$\\
$w_{1,1}$ & $w_{9,10}$ & $0$ & $0$ & $1$ & $0$ & $0$ & $1$ & $0$ & $0$ & $0$ &
$0$ & $0$ \  & $0$\\
$w_{2,1}$ & $w_{8,1}$ & $2$ & $1$ & $0$ & $2$ & $2$ & $0$ & $0$ & $0$ & $1$ &
$0$ & $0$ \  & $0$\\
$w_{2,1}$ & $w_{8,2}$ & $2$ & $0$ & $1$ & $2$ & $0$ & $0$ & $2$ & $0$ & $0$ &
$0$ & $0$ \  & $0$\\
$w_{2,1}$ & $w_{8,3}$ & $2$ & $1$ & $1$ & $0$ & $2$ & $0$ & $2$ & $0$ & $0$ &
$0$ & $1$ \  & $0$\\
$w_{2,1}$ & $w_{8,4}$ & $0$ & $0$ & $0$ & $1$ & $1$ & $0$ & $0$ & $0$ & $2$ &
$0$ & $0$ \  & $1$\\
$w_{2,1}$ & $w_{8,5}$ & $2$ & $0$ & $0$ & $1$ & $0$ & $1$ & $0$ & $0$ & $0$ &
$1$ & $0$ \  & $0$\\
$w_{2,1}$ & $w_{8,6}$ & $2$ & $2$ & $0$ & $0$ & $1$ & $1$ & $0$ & $1$ & $0$ &
$2$ & $0$ \  & $0$\\
$w_{2,1}$ & $w_{8,7}$ & $0$ & $0$ & $0$ & $2$ & $2$ & $0$ & $2$ & $0$ & $1$ &
$0$ & $1$ \  & $0$\\
$w_{2,1}$ & $w_{8,8}$ & $2$ & $0$ & $2$ & $0$ & $0$ & $2$ & $1$ & $0$ & $0$ &
$1$ & $0$ \  & $0$\\
$w_{3,1}$ & $w_{7,1}$ & $1$ & $0$ & $0$ & $1$ & $2$ & $0$ & $0$ & $0$ & $1$ &
$0$ & $0$ \  & $1$\\
$w_{3,1}$ & $w_{7,2}$ & $1$ & $0$ & $0$ & $1$ & $0$ & $0$ & $1$ & $0$ & $0$ &
$0$ & $0$ \  & $0$\\
$w_{3,1}$ & $w_{7,3}$ & $2$ & $1$ & $1$ & $3$ & $3$ & $0$ & $3$ & $0$ & $1$ &
$0$ & $1$ \  & $0$\\
$w_{3,1}$ & $w_{7,4}$ & $3$ & $2$ & $0$ & $1$ & $1$ & $1$ & $0$ & $0$ & $1$ &
$1$ & $0$ \  & $0$\\
$w_{3,1}$ & $w_{7,5}$ & $0$ & $0$ & $0$ & $1$ & $1$ & $0$ & $1$ & $0$ & $2$ &
$0$ & $1$ \  & $0$\\
$w_{3,1}$ & $w_{7,6}$ & $3$ & $0$ & $1$ & $2$ & $0$ & $1$ & $1$ & $0$ & $0$ &
$1$ & $0$ \  & $0$\\
$w_{3,1}$ & $w_{7,7}$ & $3$ & $1$ & $2$ & $0$ & $2$ & $1$ & $1$ & $1$ & $0$ &
$1$ & $0$ \  & $0$\\
$w_{3,2}$ & $w_{7,1}$ & $1$ & $1$ & $0$ & $2$ & $1$ & $0$ & $0$ & $0$ & $2$ &
$0$ & $0$ \  & $0$\\
$w_{3,2}$ & $w_{7,2}$ & $1$ & $0$ & $1$ & $1$ & $0$ & $0$ & $1$ & $0$ & $0$ &
$0$ & $0$ \  & $0$\\
$w_{3,2}$ & $w_{7,3}$ & $4$ & $1$ & $1$ & $3$ & $3$ & $0$ & $3$ & $0$ & $1$ &
$0$ & $1$ \  & $0$\\
$w_{3,2}$ & $w_{7,4}$ & $3$ & $1$ & $0$ & $2$ & $2$ & $1$ & $0$ & $1$ & $0$ &
$2$ & $0$ \  & $0$\\
$w_{3,2}$ & $w_{7,5}$ & $0$ & $0$ & $0$ & $2$ & $2$ & $0$ & $1$ & $0$ & $1$ &
$0$ & $0$ \  & $1$\\
$w_{3,2}$ & $w_{7,6}$ & $3$ & $0$ & $2$ & $1$ & $0$ & $2$ & $2$ & $0$ & $0$ &
$1$ & $0$ \  & $0$\\
$w_{3,2}$ & $w_{7,7}$ & $3$ & $2$ & $1$ & $0$ & $1$ & $2$ & $2$ & $0$ & $0$ &
$2$ & $1$ \  & $0$\\
$w_{4,1}$ & $w_{6,1}$ & $0$ & $0$ & $0$ & $1$ & $1$ & $0$ & $1$ & $0$ & $0$ &
$0$ & $0$ \  & $0$\\
$w_{4,1}$ & $w_{6,2}$ & $0$ & $0$ & $0$ & $1$ & $1$ & $0$ & $1$ & $0$ & $1$ &
$0$ & $1$ \  & $0$\\
$w_{4,1}$ & $w_{6,3}$ & $1$ & $0$ & $0$ & $0$ & $1$ & $0$ & $0$ & $0$ & $0$ &
$0$ & $0$ \  & $1$\\
$w_{4,1}$ & $w_{6,4}$ & $1$ & $0$ & $0$ & $1$ & $0$ & $0$ & $0$ & $0$ & $0$ &
$0$ & $0$ \  & $0$\\
$w_{4,1}$ & $w_{6,5}$ & $2$ & $1$ & $1$ & $1$ & $1$ & $0$ & $1$ & $0$ & $1$ &
$0$ & $0$ \  & $0$\\
$w_{4,2}$ & $w_{6,1}$ & $3$ & $1$ & $1$ & $3$ & $2$ & $0$ & $3$ & $0$ & $1$ &
$0$ & $1$ \  & $0$\\
$w_{4,2}$ & $w_{6,2}$ & $3$ & $1$ & $1$ & $4$ & $5$ & $0$ & $3$ & $0$ & $3$ &
$0$ & $1$ \  & $1$\\
$w_{4,2}$ & $w_{6,3}$ & $3$ & $2$ & $0$ & $2$ & $2$ & $1$ & $0$ & $0$ & $2$ &
$1$ & $0$ \  & $0$\\
$w_{4,2}$ & $w_{6,4}$ & $3$ & $0$ & $1$ & $2$ & $0$ & $1$ & $2$ & $0$ & $0$ &
$1$ & $0$ \  & $0$\\
$w_{4,2}$ & $w_{6,5}$ & $9$ & $3$ & $3$ & $5$ & $5$ & $3$ & $4$ & $1$ & $1$ &
$3$ & $1$ \  & $0$\\
$w_{4,3}$ & $w_{6,1}$ & $2$ & $0$ & $1$ & $1$ & $1$ & $0$ & $1$ & $0$ & $0$ &
$0$ & $0$ \  & $0$\\
$w_{4,3}$ & $w_{6,2}$ & $2$ & $1$ & $0$ & $3$ & $1$ & $0$ & $1$ & $0$ & $1$ &
$0$ & $0$ \  & $0$\\
$w_{4,3}$ & $w_{6,3}$ & $1$ & $0$ & $0$ & $2$ & $1$ & $0$ & $0$ & $1$ & $0$ &
$1$ & $0$ \  & $0$\\
$w_{4,3}$ & $w_{6,4}$ & $1$ & $0$ & $2$ & $0$ & $0$ & $1$ & $1$ & $0$ & $0$ &
$0$ & $0$ \  & $0$\\
$w_{4,3}$ & $w_{6,5}$ & $4$ & $1$ & $1$ & $1$ & $1$ & $2$ & $3$ & $0$ & $0$ &
$2$ & $1$ \  & $0$\\
$w_{5,1}$ & $w_{5,1}$ & $0$ & $0$ & $0$ & $2$ & $2$ & $0$ & $2$ & $0$ & $1$ &
$0$ & $1$ \  & $0$\\
$w_{5,1}$ & $w_{5,2}$ & $3$ & $1$ & $1$ & $3$ & $2$ & $0$ & $2$ & $0$ & $1$ &
$0$ & $0$ \  & $0$\\
$w_{5,1}$ & $w_{5,3}$ & $3$ & $1$ & $1$ & $2$ & $3$ & $0$ & $2$ & $0$ & $2$ &
$0$ & $1$ \  & $1$\\
$w_{5,1}$ & $w_{5,4}$ & $1$ & $0$ & $0$ & $1$ & $0$ & $0$ & $0$ & $0$ & $0$ &
$0$ & $0$ \  & $0$\\
$w_{5,2}$ & $w_{5,2}$ & $6$ & $1$ & $3$ & $2$ & $2$ & $2$ & $4$ & $0$ & $0$ &
$2$ & $1$ \  & $0$\\
$w_{5,2}$ & $w_{5,3}$ & $6$ & $2$ & $1$ & $5$ & $3$ & $2$ & $3$ & $1$ & $1$ &
$2$ & $1$ \  & $0$\\
$w_{5,2}$ & $w_{5,4}$ & $1$ & $0$ & $1$ & $0$ & $0$ & $1$ & $1$ & $0$ & $0$ &
$0$ & $0$ \  & $0$\\
$w_{5,3}$ & $w_{5,3}$ & $6$ & $3$ & $2$ & $4$ & $6$ & $2$ & $2$ & $0$ & $3$ &
$2$ & $0$ \  & $0$\\
$w_{5,3}$ & $w_{5,4}$ & $1$ & $0$ & $0$ & $1$ & $0$ & $0$ & $1$ & $0$ & $0$ &
$1$ & $0$ \  & $0$\\
$w_{5,4}$ & $w_{5,4}$ & $0$ & $0$ & $1$ & $0$ & $0$ & $0$ & $0$ & $0$ & $0$ &
$0$ & $0$ \  & $0$\\\hline
\end{tabular}
}

\begin{center}
\end{center}

\newpage\renewcommand{\arraystretch}{1}

\begin{center}
{ }
\end{center}

\textbf{Table B}$_{8}${\footnotesize . }\textsl{L-R coefficients for }%
$E_{8}/S^{1}\cdot SU(8)$

{\footnotesize \setlength{\tabcolsep}{1mm}
\begin{tabular}
[c]{|r|r|r|r|r|r|r|r|r|r|r|r|r|r|r|}\hline
$u$ & $v$ & \multicolumn{13}{c|}{$w\in W^{9}(8)$}\\\cline{3-15}
&  & $w_{9,1}$ & $w_{9,2}$ & $w_{9,3}$ & $w_{9,4}$ & $w_{9,5}$ & $w_{9,6}$ &
$w_{9,7}$ & $w_{9,8}$ & $w_{9,9}$ & $w_{9,10}$ & $w_{9,11}$ & $w_{9,12}$ &
\ $w_{9,13}$\\\hline
$w_{1,1}$ & $w_{8,1}$ & $1$ & $1$ & $0$ & $0$ & $1$ & $0$ & $0$ & $0$ & $0$ &
$0$ & $0$ \  & $0$ & $0$\\
$w_{1,1}$ & $w_{8,2}$ & $0$ & $0$ & $1$ & $0$ & $0$ & $0$ & $0$ & $0$ & $0$ &
$0$ & $0$ \  & $0$ & $0$\\
$w_{1,1}$ & $w_{8,3}$ & $1$ & $0$ & $1$ & $1$ & $0$ & $0$ & $0$ & $1$ & $0$ &
$0$ & $0$ \  & $0$ & $0$\\
$w_{1,1}$ & $w_{8,4}$ & $0$ & $1$ & $0$ & $1$ & $0$ & $0$ & $0$ & $0$ & $1$ &
$0$ & $0$ \  & $0$ & $0$\\
$w_{1,1}$ & $w_{8,5}$ & $0$ & $0$ & $0$ & $0$ & $1$ & $0$ & $0$ & $0$ & $0$ &
$1$ & $0$ \  & $0$ & $0$\\
$w_{1,1}$ & $w_{8,6}$ & $1$ & $0$ & $0$ & $0$ & $0$ & $1$ & $1$ & $0$ & $0$ &
$0$ & $0$ \  & $0$ & $0$\\
$w_{1,1}$ & $w_{8,7}$ & $0$ & $1$ & $0$ & $0$ & $0$ & $0$ & $1$ & $0$ & $0$ &
$0$ & $1$ \  & $0$ & $0$\\
$w_{1,1}$ & $w_{8,8}$ & $0$ & $0$ & $0$ & $0$ & $1$ & $0$ & $0$ & $1$ & $1$ &
$0$ & $0$ \  & $0$ & $0$\\
$w_{1,1}$ & $w_{8,9}$ & $0$ & $0$ & $1$ & $0$ & $0$ & $1$ & $0$ & $0$ & $0$ &
$0$ & $0$ \  & $1$ & $0$\\
$w_{1,1}$ & $w_{8,10}$ & $0$ & $0$ & $0$ & $1$ & $0$ & $0$ & $1$ & $0$ & $0$ &
$0$ & $0$ \  & $1$ & $1$\\
$w_{2,1}$ & $w_{7,1}$ & $1$ & $1$ & $0$ & $0$ & $2$ & $0$ & $0$ & $0$ & $0$ &
$1$ & $0$ \  & $0$ & $0$\\
$w_{2,1}$ & $w_{7,2}$ & $1$ & $0$ & $2$ & $1$ & $0$ & $0$ & $0$ & $1$ & $0$ &
$0$ & $0$ \  & $0$ & $0$\\
$w_{2,1}$ & $w_{7,3}$ & $2$ & $2$ & $1$ & $2$ & $2$ & $0$ & $0$ & $2$ & $2$ &
$0$ & $0$ \  & $0$ & $0$\\
$w_{2,1}$ & $w_{7,4}$ & $2$ & $2$ & $0$ & $0$ & $1$ & $1$ & $2$ & $0$ & $0$ &
$0$ & $1$ \  & $0$ & $0$\\
$w_{2,1}$ & $w_{7,5}$ & $0$ & $0$ & $0$ & $0$ & $2$ & $0$ & $0$ & $1$ & $1$ &
$1$ & $0$ \  & $0$ & $0$\\
$w_{2,1}$ & $w_{7,6}$ & $0$ & $0$ & $2$ & $0$ & $0$ & $1$ & $0$ & $0$ & $0$ &
$0$ & $0$ \  & $1$ & $0$\\
$w_{2,1}$ & $w_{7,7}$ & $2$ & $0$ & $2$ & $2$ & $0$ & $2$ & $2$ & $1$ & $0$ &
$0$ & $0$ \  & $2$ & $1$\\
$w_{2,1}$ & $w_{7,8}$ & $0$ & $2$ & $0$ & $2$ & $0$ & $0$ & $2$ & $0$ & $1$ &
$0$ & $1$ \  & $1$ & $1$\\
$w_{3,1}$ & $w_{6,1}$ & $1$ & $1$ & $1$ & $1$ & $1$ & $0$ & $0$ & $2$ & $1$ &
$0$ & $0$ \  & $0$ & $0$\\
$w_{3,1}$ & $w_{6,2}$ & $1$ & $1$ & $0$ & $1$ & $3$ & $0$ & $0$ & $1$ & $2$ &
$1$ & $0$ \  & $0$ & $0$\\
$w_{3,1}$ & $w_{6,3}$ & $1$ & $2$ & $0$ & $0$ & $1$ & $0$ & $1$ & $0$ & $0$ &
$1$ & $0$ \  & $0$ & $0$\\
$w_{3,1}$ & $w_{6,4}$ & $2$ & $0$ & $3$ & $1$ & $0$ & $1$ & $1$ & $1$ & $0$ &
$0$ & $0$ \  & $1$ & $0$\\
$w_{3,1}$ & $w_{6,5}$ & $3$ & $3$ & $1$ & $3$ & $2$ & $1$ & $2$ & $1$ & $1$ &
$0$ & $1$ \  & $1$ & $1$\\
$w_{3,1}$ & $w_{6,6}$ & $0$ & $0$ & $1$ & $0$ & $0$ & $1$ & $0$ & $0$ & $0$ &
$0$ & $0$ \  & $0$ & $0$\\
$w_{3,2}$ & $w_{6,1}$ & $2$ & $1$ & $2$ & $2$ & $1$ & $0$ & $0$ & $1$ & $1$ &
$0$ & $0$ \  & $0$ & $0$\\
$w_{3,2}$ & $w_{6,2}$ & $2$ & $2$ & $1$ & $1$ & $3$ & $0$ & $0$ & $2$ & $1$ &
$1$ & $0$ \  & $0$ & $0$\\
$w_{3,2}$ & $w_{6,3}$ & $2$ & $1$ & $0$ & $0$ & $2$ & $1$ & $1$ & $0$ & $0$ &
$0$ & $1$ \  & $0$ & $0$\\
$w_{3,2}$ & $w_{6,4}$ & $1$ & $0$ & $3$ & $2$ & $0$ & $2$ & $1$ & $1$ & $0$ &
$0$ & $0$ \  & $2$ & $1$\\
$w_{3,2}$ & $w_{6,5}$ & $3$ & $3$ & $2$ & $3$ & $1$ & $2$ & $4$ & $2$ & $2$ &
$0$ & $1$ \  & $2$ & $1$\\
$w_{3,2}$ & $w_{6,6}$ & $0$ & $0$ & $1$ & $0$ & $0$ & $0$ & $0$ & $0$ & $0$ &
$0$ & $0$ \  & $1$ & $0$\\
$w_{4,1}$ & $w_{5,1}$ & $0$ & $0$ & $0$ & $0$ & $1$ & $0$ & $0$ & $1$ & $1$ &
$0$ & $0$ \  & $0$ & $0$\\
$w_{4,1}$ & $w_{5,2}$ & $1$ & $1$ & $1$ & $1$ & $1$ & $0$ & $0$ & $1$ & $0$ &
$0$ & $0$ \  & $0$ & $0$\\
$w_{4,1}$ & $w_{5,3}$ & $1$ & $1$ & $0$ & $1$ & $1$ & $0$ & $0$ & $0$ & $1$ &
$1$ & $0$ \  & $0$ & $0$\\
$w_{4,1}$ & $w_{5,4}$ & $1$ & $0$ & $1$ & $0$ & $0$ & $0$ & $0$ & $0$ & $0$ &
$0$ & $0$ \  & $0$ & $0$\\
$w_{4,2}$ & $w_{5,1}$ & $2$ & $2$ & $1$ & $2$ & $3$ & $0$ & $0$ & $2$ & $2$ &
$1$ & $0$ \  & $0$ & $0$\\
$w_{4,2}$ & $w_{5,2}$ & $5$ & $3$ & $4$ & $4$ & $2$ & $2$ & $3$ & $3$ & $2$ &
$0$ & $1$ \  & $2$ & $1$\\
$w_{4,2}$ & $w_{5,3}$ & $4$ & $5$ & $1$ & $3$ & $5$ & $1$ & $3$ & $2$ & $2$ &
$1$ & $1$ \  & $1$ & $1$\\
$w_{4,2}$ & $w_{5,4}$ & $1$ & $0$ & $3$ & $1$ & $0$ & $2$ & $1$ & $1$ & $0$ &
$0$ & $0$ \  & $1$ & $0$\\
$w_{4,3}$ & $w_{5,1}$ & $2$ & $1$ & $2$ & $1$ & $1$ & $0$ & $0$ & $1$ & $0$ &
$0$ & $0$ \  & $0$ & $0$\\
$w_{4,3}$ & $w_{5,2}$ & $1$ & $1$ & $3$ & $3$ & $0$ & $2$ & $2$ & $1$ & $1$ &
$0$ & $0$ \  & $2$ & $1$\\
$w_{4,3}$ & $w_{5,3}$ & $3$ & $1$ & $2$ & $1$ & $1$ & $2$ & $2$ & $2$ & $1$ &
$0$ & $1$ \  & $1$ & $0$\\
$w_{4,3}$ & $w_{5,4}$ & $0$ & $0$ & $1$ & $1$ & $0$ & $0$ & $0$ & $0$ & $0$ &
$0$ & $0$ \  & $2$ & $1$\\\hline
\end{tabular}
}

\begin{center}
\end{center}

\newpage\renewcommand{\arraystretch}{0.8}

\begin{center}
{ \setlength{\tabcolsep}{0.2mm}
{\tiny
\begin{tabular}
[c]{|r|r|r|r|r|r|r|r|r|r|r|r|r|r|r|r|r|r|r|}\hline
$u$ & $v$ & \multicolumn{17}{c|}{$w\in W^{10}(8)$}\\\cline{3-19}
&  & $w_{10,1}$ & $w_{10,2}$ & $w_{10,3}$ & $w_{10,4}$ & $w_{10,5}$ &
$w_{10,6}$ & $w_{10,7}$ & \ $w_{10,8}$ & $w_{10,9}$ & $w_{10,10}$ &
$w_{10,11}$ & \ $w_{10,12}$ & $w_{10,13}$ & $w_{10,14}$ & $w_{10,15}$ &
$w_{10,16}$ & $w_{10,17}$\\\hline
$w_{1,1}$ & $w_{9,1}$ & $1$ & $1$ & $0$ & $0$ & $0$ & $1$ & $0$ & $0$ & $0$ &
$0$ & $0$ \  & $0$ & $0$ & $0$ & $0$ & $0$ & $0$\\
$w_{1,1}$ & $w_{9,2}$ & $0$ & $1$ & $1$ & $0$ & $0$ & $0$ & $1$ & $0$ & $0$ &
$0$ & $0$ \  & $0$ & $0$ & $0$ & $0$ & $0$ & $0$\\
$w_{1,1}$ & $w_{9,3}$ & $1$ & $0$ & $0$ & $1$ & $0$ & $0$ & $0$ & $0$ & $0$ &
$1$ & $0$ \  & $0$ & $0$ & $0$ & $0$ & $0$ & $0$\\
$w_{1,1}$ & $w_{9,4}$ & $0$ & $1$ & $0$ & $1$ & $1$ & $0$ & $0$ & $0$ & $0$ &
$0$ & $1$ \  & $0$ & $0$ & $0$ & $0$ & $0$ & $0$\\
$w_{1,1}$ & $w_{9,5}$ & $0$ & $0$ & $0$ & $0$ & $0$ & $1$ & $1$ & $0$ & $0$ &
$0$ & $0$ \  & $0$ & $1$ & $0$ & $0$ & $0$ & $0$\\
$w_{1,1}$ & $w_{9,6}$ & $1$ & $0$ & $0$ & $0$ & $0$ & $0$ & $0$ & $1$ & $0$ &
$0$ & $0$ \  & $0$ & $0$ & $0$ & $0$ & $0$ & $0$\\
$w_{1,1}$ & $w_{9,7}$ & $0$ & $1$ & $0$ & $0$ & $0$ & $0$ & $0$ & $1$ & $1$ &
$0$ & $0$ \  & $0$ & $0$ & $1$ & $0$ & $0$ & $0$\\
$w_{1,1}$ & $w_{9,8}$ & $0$ & $0$ & $0$ & $0$ & $0$ & $1$ & $0$ & $0$ & $0$ &
$1$ & $1$ \  & $0$ & $0$ & $0$ & $0$ & $0$ & $0$\\
$w_{1,1}$ & $w_{9,9}$ & $0$ & $0$ & $0$ & $0$ & $0$ & $0$ & $1$ & $0$ & $0$ &
$0$ & $1$ \  & $0$ & $0$ & $0$ & $1$ & $0$ & $0$\\
$w_{1,1}$ & $w_{9,10}$ & $0$ & $0$ & $0$ & $0$ & $0$ & $0$ & $0$ & $0$ & $0$ &
$0$ & $0$ \  & $0$ & $1$ & $0$ & $0$ & $1$ & $0$\\
$w_{1,1}$ & $w_{9,11}$ & $0$ & $0$ & $1$ & $0$ & $0$ & $0$ & $0$ & $0$ & $0$ &
$0$ & $0$ \  & $1$ & $0$ & $1$ & $0$ & $0$ & $0$\\
$w_{1,1}$ & $w_{9,12}$ & $0$ & $0$ & $0$ & $1$ & $0$ & $0$ & $0$ & $1$ & $0$ &
$0$ & $0$ \  & $0$ & $0$ & $0$ & $0$ & $0$ & $1$\\
$w_{1,1}$ & $w_{9,13}$ & $0$ & $0$ & $0$ & $0$ & $1$ & $0$ & $0$ & $0$ & $1$ &
$0$ & $0$ \  & $0$ & $0$ & $0$ & $0$ & $0$ & $1$\\
$w_{2,1}$ & $w_{8,1}$ & $1$ & $2$ & $1$ & $0$ & $0$ & $2$ & $2$ & $0$ & $0$ &
$0$ & $0$ \  & $0$ & $1$ & $0$ & $0$ & $0$ & $0$\\
$w_{2,1}$ & $w_{8,2}$ & $1$ & $0$ & $0$ & $1$ & $0$ & $0$ & $0$ & $0$ & $0$ &
$1$ & $0$ \  & $0$ & $0$ & $0$ & $0$ & $0$ & $0$\\
$w_{2,1}$ & $w_{8,3}$ & $2$ & $2$ & $0$ & $2$ & $1$ & $2$ & $0$ & $0$ & $0$ &
$2$ & $2$ \  & $0$ & $0$ & $0$ & $0$ & $0$ & $0$\\
$w_{2,1}$ & $w_{8,4}$ & $0$ & $2$ & $1$ & $1$ & $1$ & $0$ & $2$ & $0$ & $0$ &
$0$ & $2$ \  & $0$ & $0$ & $0$ & $1$ & $0$ & $0$\\
$w_{2,1}$ & $w_{8,5}$ & $0$ & $0$ & $0$ & $0$ & $0$ & $1$ & $1$ & $0$ & $0$ &
$0$ & $0$ \  & $0$ & $2$ & $0$ & $0$ & $1$ & $0$\\
$w_{2,1}$ & $w_{8,6}$ & $2$ & $2$ & $0$ & $0$ & $0$ & $1$ & $0$ & $2$ & $1$ &
$0$ & $0$ \  & $0$ & $0$ & $1$ & $0$ & $0$ & $0$\\
$w_{2,1}$ & $w_{8,7}$ & $0$ & $2$ & $2$ & $0$ & $0$ & $0$ & $1$ & $1$ & $1$ &
$0$ & $0$ \  & $1$ & $0$ & $2$ & $0$ & $0$ & $0$\\
$w_{2,1}$ & $w_{8,8}$ & $0$ & $0$ & $0$ & $0$ & $0$ & $2$ & $2$ & $0$ & $0$ &
$1$ & $2$ \  & $0$ & $1$ & $0$ & $1$ & $0$ & $0$\\
$w_{2,1}$ & $w_{8,9}$ & $2$ & $0$ & $0$ & $2$ & $0$ & $0$ & $0$ & $2$ & $0$ &
$1$ & $0$ \  & $0$ & $0$ & $0$ & $0$ & $0$ & $1$\\
$w_{2,1}$ & $w_{8,10}$ & $0$ & $2$ & $0$ & $2$ & $2$ & $0$ & $0$ & $2$ & $2$ &
$0$ & $1$ \  & $0$ & $0$ & $1$ & $0$ & $0$ & $2$\\
$w_{3,1}$ & $w_{7,1}$ & $0$ & $1$ & $0$ & $0$ & $0$ & $1$ & $2$ & $0$ & $0$ &
$0$ & $0$ \  & $0$ & $1$ & $0$ & $0$ & $1$ & $0$\\
$w_{3,1}$ & $w_{7,2}$ & $1$ & $1$ & $0$ & $1$ & $0$ & $1$ & $0$ & $0$ & $0$ &
$2$ & $1$ \  & $0$ & $0$ & $0$ & $0$ & $0$ & $0$\\
$w_{3,1}$ & $w_{7,3}$ & $1$ & $2$ & $1$ & $1$ & $1$ & $3$ & $3$ & $0$ & $0$ &
$1$ & $3$ \  & $0$ & $1$ & $0$ & $1$ & $0$ & $0$\\
$w_{3,1}$ & $w_{7,4}$ & $1$ & $3$ & $2$ & $0$ & $0$ & $1$ & $1$ & $1$ & $1$ &
$0$ & $0$ \  & $0$ & $1$ & $1$ & $0$ & $0$ & $0$\\
$w_{3,1}$ & $w_{7,5}$ & $0$ & $0$ & $0$ & $0$ & $0$ & $1$ & $1$ & $0$ & $0$ &
$0$ & $1$ \  & $0$ & $2$ & $0$ & $1$ & $0$ & $0$\\
$w_{3,1}$ & $w_{7,6}$ & $2$ & $0$ & $0$ & $1$ & $0$ & $0$ & $0$ & $1$ & $0$ &
$1$ & $0$ \  & $0$ & $0$ & $0$ & $0$ & $0$ & $0$\\
$w_{3,1}$ & $w_{7,7}$ & $3$ & $3$ & $0$ & $3$ & $1$ & $2$ & $0$ & $2$ & $1$ &
$1$ & $1$ \  & $0$ & $0$ & $1$ & $0$ & $0$ & $1$\\
$w_{3,1}$ & $w_{7,8}$ & $0$ & $3$ & $1$ & $1$ & $2$ & $0$ & $2$ & $1$ & $1$ &
$0$ & $1$ \  & $1$ & $0$ & $1$ & $0$ & $0$ & $1$\\
$w_{3,2}$ & $w_{7,1}$ & $1$ & $1$ & $1$ & $0$ & $0$ & $2$ & $1$ & $0$ & $0$ &
$0$ & $0$ \  & $0$ & $2$ & $0$ & $0$ & $0$ & $0$\\
$w_{3,2}$ & $w_{7,2}$ & $2$ & $1$ & $0$ & $2$ & $1$ & $1$ & $0$ & $0$ & $0$ &
$1$ & $1$ \  & $0$ & $0$ & $0$ & $0$ & $0$ & $0$\\
$w_{3,2}$ & $w_{7,3}$ & $2$ & $4$ & $1$ & $2$ & $1$ & $3$ & $3$ & $0$ & $0$ &
$2$ & $3$ \  & $0$ & $1$ & $0$ & $1$ & $0$ & $0$\\
$w_{3,2}$ & $w_{7,4}$ & $2$ & $3$ & $1$ & $0$ & $0$ & $2$ & $2$ & $2$ & $1$ &
$0$ & $0$ \  & $1$ & $0$ & $2$ & $0$ & $0$ & $0$\\
$w_{3,2}$ & $w_{7,5}$ & $0$ & $0$ & $0$ & $0$ & $0$ & $2$ & $2$ & $0$ & $0$ &
$1$ & $1$ \  & $0$ & $1$ & $0$ & $0$ & $1$ & $0$\\
$w_{3,2}$ & $w_{7,6}$ & $1$ & $0$ & $0$ & $2$ & $0$ & $0$ & $0$ & $1$ & $0$ &
$1$ & $0$ \  & $0$ & $0$ & $0$ & $0$ & $0$ & $1$\\
$w_{3,2}$ & $w_{7,7}$ & $3$ & $3$ & $0$ & $3$ & $2$ & $1$ & $0$ & $4$ & $2$ &
$2$ & $2$ \  & $0$ & $0$ & $1$ & $0$ & $0$ & $2$\\
$w_{3,2}$ & $w_{7,8}$ & $0$ & $3$ & $2$ & $2$ & $1$ & $0$ & $1$ & $2$ & $2$ &
$0$ & $2$ \  & $0$ & $0$ & $2$ & $1$ & $0$ & $1$\\
$w_{4,1}$ & $w_{6,1}$ & $0$ & $0$ & $0$ & $0$ & $0$ & $1$ & $1$ & $0$ & $0$ &
$1$ & $1$ \  & $0$ & $0$ & $0$ & $0$ & $0$ & $0$\\
$w_{4,1}$ & $w_{6,2}$ & $0$ & $0$ & $0$ & $0$ & $0$ & $1$ & $1$ & $0$ & $0$ &
$0$ & $1$ \  & $0$ & $1$ & $0$ & $1$ & $0$ & $0$\\
$w_{4,1}$ & $w_{6,3}$ & $0$ & $1$ & $0$ & $0$ & $0$ & $0$ & $1$ & $0$ & $0$ &
$0$ & $0$ \  & $0$ & $0$ & $0$ & $0$ & $1$ & $0$\\
$w_{4,1}$ & $w_{6,4}$ & $1$ & $1$ & $0$ & $1$ & $0$ & $1$ & $0$ & $0$ & $0$ &
$1$ & $0$ \  & $0$ & $0$ & $0$ & $0$ & $0$ & $0$\\
$w_{4,1}$ & $w_{6,5}$ & $1$ & $2$ & $1$ & $1$ & $1$ & $1$ & $1$ & $0$ & $0$ &
$0$ & $1$ \  & $0$ & $1$ & $0$ & $0$ & $0$ & $0$\\
$w_{4,1}$ & $w_{6,6}$ & $1$ & $0$ & $0$ & $0$ & $0$ & $0$ & $0$ & $0$ & $0$ &
$0$ & $0$ \  & $0$ & $0$ & $0$ & $0$ & $0$ & $0$\\
$w_{4,2}$ & $w_{6,1}$ & $2$ & $3$ & $1$ & $2$ & $1$ & $3$ & $2$ & $0$ & $0$ &
$2$ & $3$ \  & $0$ & $1$ & $0$ & $1$ & $0$ & $0$\\
$w_{4,2}$ & $w_{6,2}$ & $1$ & $3$ & $1$ & $1$ & $1$ & $4$ & $5$ & $0$ & $0$ &
$1$ & $3$ \  & $0$ & $3$ & $0$ & $1$ & $1$ & $0$\\
$w_{4,2}$ & $w_{6,3}$ & $1$ & $3$ & $2$ & $0$ & $0$ & $2$ & $2$ & $1$ & $1$ &
$0$ & $0$ \  & $0$ & $2$ & $1$ & $0$ & $0$ & $0$\\
$w_{4,2}$ & $w_{6,4}$ & $5$ & $3$ & $0$ & $4$ & $1$ & $2$ & $0$ & $3$ & $1$ &
$3$ & $2$ \  & $0$ & $0$ & $1$ & $0$ & $0$ & $1$\\
$w_{4,2}$ & $w_{6,5}$ & $4$ & $9$ & $3$ & $4$ & $3$ & $5$ & $5$ & $4$ & $3$ &
$2$ & $4$ \  & $1$ & $1$ & $3$ & $1$ & $0$ & $2$\\
$w_{4,2}$ & $w_{6,6}$ & $1$ & $0$ & $0$ & $1$ & $0$ & $0$ & $0$ & $1$ & $0$ &
$1$ & $0$ \  & $0$ & $0$ & $0$ & $0$ & $0$ & $0$\\
$w_{4,3}$ & $w_{6,1}$ & $2$ & $2$ & $0$ & $2$ & $1$ & $1$ & $1$ & $0$ & $0$ &
$1$ & $1$ \  & $0$ & $0$ & $0$ & $0$ & $0$ & $0$\\
$w_{4,3}$ & $w_{6,2}$ & $2$ & $2$ & $1$ & $1$ & $0$ & $3$ & $1$ & $0$ & $0$ &
$2$ & $1$ \  & $0$ & $1$ & $0$ & $0$ & $0$ & $0$\\
$w_{4,3}$ & $w_{6,3}$ & $2$ & $1$ & $0$ & $0$ & $0$ & $2$ & $1$ & $1$ & $0$ &
$0$ & $0$ \  & $1$ & $0$ & $1$ & $0$ & $0$ & $0$\\
$w_{4,3}$ & $w_{6,4}$ & $1$ & $1$ & $0$ & $3$ & $2$ & $0$ & $0$ & $2$ & $1$ &
$1$ & $1$ \  & $0$ & $0$ & $0$ & $0$ & $0$ & $2$\\
$w_{4,3}$ & $w_{6,5}$ & $3$ & $4$ & $1$ & $3$ & $1$ & $1$ & $1$ & $4$ & $2$ &
$2$ & $3$ \  & $0$ & $0$ & $2$ & $1$ & $0$ & $1$\\
$w_{4,3}$ & $w_{6,6}$ & $0$ & $0$ & $0$ & $1$ & $0$ & $0$ & $0$ & $0$ & $0$ &
$0$ & $0$ \  & $0$ & $0$ & $0$ & $0$ & $0$ & $1$\\
$w_{5,1}$ & $w_{5,1}$ & $0$ & $0$ & $0$ & $0$ & $0$ & $2$ & $2$ & $0$ & $0$ &
$1$ & $2$ \  & $0$ & $1$ & $0$ & $1$ & $0$ & $0$\\
$w_{5,1}$ & $w_{5,2}$ & $2$ & $3$ & $1$ & $2$ & $1$ & $3$ & $2$ & $0$ & $0$ &
$2$ & $2$ \  & $0$ & $1$ & $0$ & $0$ & $0$ & $0$\\
$w_{5,1}$ & $w_{5,3}$ & $1$ & $3$ & $1$ & $1$ & $1$ & $2$ & $3$ & $0$ & $0$ &
$0$ & $2$ \  & $0$ & $2$ & $0$ & $1$ & $1$ & $0$\\
$w_{5,1}$ & $w_{5,4}$ & $2$ & $1$ & $0$ & $1$ & $0$ & $1$ & $0$ & $0$ & $0$ &
$1$ & $0$ \  & $0$ & $0$ & $0$ & $0$ & $0$ & $0$\\
$w_{5,2}$ & $w_{5,2}$ & $5$ & $6$ & $1$ & $5$ & $3$ & $2$ & $2$ & $4$ & $2$ &
$3$ & $4$ \  & $0$ & $0$ & $2$ & $1$ & $0$ & $2$\\
$w_{5,2}$ & $w_{5,3}$ & $4$ & $6$ & $2$ & $3$ & $1$ & $5$ & $3$ & $3$ & $2$ &
$2$ & $3$ \  & $1$ & $1$ & $2$ & $1$ & $0$ & $1$\\
$w_{5,2}$ & $w_{5,4}$ & $1$ & $1$ & $0$ & $3$ & $1$ & $0$ & $0$ & $2$ & $1$ &
$1$ & $1$ \  & $0$ & $0$ & $0$ & $0$ & $0$ & $1$\\
$w_{5,3}$ & $w_{5,3}$ & $1$ & $6$ & $3$ & $1$ & $2$ & $4$ & $6$ & $2$ & $2$ &
$1$ & $2$ \  & $0$ & $3$ & $2$ & $0$ & $0$ & $1$\\
$w_{5,3}$ & $w_{5,4}$ & $3$ & $1$ & $0$ & $1$ & $0$ & $1$ & $0$ & $2$ & $0$ &
$2$ & $1$ \  & $0$ & $0$ & $1$ & $0$ & $0$ & $0$\\
$w_{5,4}$ & $w_{5,4}$ & $0$ & $0$ & $0$ & $1$ & $1$ & $0$ & $0$ & $0$ & $0$ &
$0$ & $0$ \  & $0$ & $0$ & $0$ & $0$ & $0$ & $2$\\\hline
\end{tabular}
}}
\end{center}

The computations were carried out by using Mathematicae on a PC. PIV667. Ram
128. Win98. In general, the running time of the program depends on

\begin{quote}
(1) the order of the coset $\overline{W}$;

(2) the number of non-zero entries in the Cartan matrix of $G$.
\end{quote}

\noindent More precisely, to obtain the results in Table A$_{n}$, the times
consumed (in seconds) are \renewcommand{\arraystretch}{1.2}

\begin{center}%
\begin{tabular}
[c]{|c|c|c|c|}\hline
$n$ & 6 & 7 & 8\\\hline
{time} & 1 & 1 & 2\\\hline
\end{tabular}
.
\end{center}

\noindent The running times for computing all $a_{u,v}^{w}$ with $l(w)\leq10$
are respectively

\begin{center}%
\begin{tabular}
[c]{|c|c|c|c|}\hline
$n$ & 6 & 7 & 8\\\hline
{time} & 38 & 115 & 159\\\hline
\end{tabular}
.

\bigskip

\textbf{References}
\end{center}

[B] P. Baum, On the cohomology of homogeneous spaces, Topology 7(1968), 15-38.

[Be] N. Bergeron, A combinatorial construction of the Schubert polynomials, J.
Combin. Theory, Ser.A, 60(1992), 168-182.

[BGG] I. N. Bernstein, I. M. Gel'fand and S. I. Gel'fand, Schubert cells and
cohomology of the spaces G/P, Russian Math. Surveys 28 (1973), 1-26.

[Bi] S. Billey, Kostant polynomials and the cohomology ring for G/B, Duke J.
Math. 96, No.1(1999), 205-224.

[BH] S. Billey and M. Haiman, Schubert polynomials for the classical groups,
J. AMS, 8 (no. 2)(1995), 443-482.

[BHi] A. Borel and F. Hirzebruch, Characteristic classes and homogeneous
spaces (I), Amer. J. Math. 80, 1958, 458--538.

[BJS] S. Billey, W. Jockush and S. Stanley, Some combinatorial properties of
Schubert polynomials, J. Algebraic Combin., 2 (no. 4)(1993), 345-375.

[Bo$_{1}$] A. Borel, Sur la cohomologie des espaces fibr\'{e}s principaux et
des espaces homogenes de groupes de Lie compacts, Ann. Math. 57(1953), 115-207.

[Bo$_{2}$] A. Borel, Topics in the homology theory of fiber bundles, Berlin,
Springer, 1967.

[BS] R. Bott and H. Samelson, Application of the theory of Morse to symmetric
spaces, Amer. J. Math., Vol. LXXX, no. 4 (1958), 964-1029.

[BS$_{1}$] N. Bergeron and F. Sottile, A Pieri-type formula for isotropic flag
manifolds, Trans. Amer. Math. Soc. 354 (7)(2002), 2659--2705.

[BS$_{2}$] N. Bergeron and F. Sottile, Skew Schubert functions and the Pieri
formula for flag manifolds, Trans. Amer. Math. Soc. 354 (2) (2002), 651--673.

[BS$_{3}$] N. Bergeron and F. Sottile, Schubert polynomials, the Bruhat order,
and the geometry of flag manifolds, Duke Math. J. 95 (2)(1998), 373--423.

[Ch$_{1}$] C. Chevalley, La th\'{e}orie des groupes alg\'{e}briques, Proc.
1958 ICM, Cambridge Univ. Press, 1960, 53-68.

[Ch$_{2}$] C. Chevalley, Sur les D\'{e}compositions Celluaires des Espaces
G/B, in Algebraic groups and their generalizations: Classical methods, W.
Haboush ed. Proc. Symp. in Pure Math. 56 (part 1) (1994), 1-26.

[D] M. Demazure, D\'{e}singularization des vari\'{e}t\'{e}s de Schubert
g\'{e}n\'{e}ralis\'{e}es, Ann. Sci. \'{E}cole. Norm. Sup. (4) 7(1974), 53-88.

[Du$_{1}$] H. Duan, Self-maps of the Grassmannian of complex structures.
Compositio Math. 132 (2002), no. 2, 159--175.

[Du$_{2}$] H. Duan, The degree of a Schubert variety, Adv. Math., 180(2003),112-133.

[Du$_{3}$] H. Duan, Multiplicative rule of Schubert classes, Invent. Math. 159
(2005), no. 2, 407-436.

[Du$_{4}$] H. Duan, Multiplicative rule in the Grothendieck cohomology of a
flag variety, arXiv: math.AG/0411588

[DZ] H. Duan and Xuezhi Zhao, A unified formula for Steenrod operations on
flag manifolds, arXiv: math.AT/0306250.

[DZZ] H. Duan, Xu-an Zhao and Xuezhi Zhao, The Cartan matrix and enumerative
calculus, J. Symbolic comput.,\textsl{ }38(2004), 1119-1144.

[E] C. Ehresmann, Sur la topologie de certains espaces homogenes, Ann. of
Math. (2) 35 (1934), 396--443.

[FK] S. Fomin and A. Kirillov, Combinatorial B$_{n}$-analogs of Schubert
polynomials, Trans. AMS 348(1996), 3591-3620.

[FS] S. Fomin and R. Stanley, Schubert polynomials and nilCoxeter algebra,
Adv. Math., 103(1994), 196-207.

[Fu] W. Fulton, Universal Schubert polynomials, Duke Math. J. 96(no. 3)(1999), 575-594.

[H] M. Hoffman, Endomorphisms of the cohomology of complex Grassmannians.
Trans. Amer. Math. Soc. 281 (1984), no. 2, 745--760.

[Han] H.C. Hansen, On cycles in flag manifolds, Math. Scand. 33 (1973), 269-274.

[HB] H. Hiller and B. Boe, Pieri formula for $SO_{2n+1}/U_{n}$ and
$Sp_{n}/U_{n}$, Adv. in Math. 62 (1)(1986), 49--67.

[Hu] J. E. Humphreys, Introduction to Lie algebras and representation theory,
Graduated Texts in Math. 9, Springer-Verlag New York, 1972.

[IM] A.Iliev and L. Manivel, The Chow ring of the Cayley plane. Compos. Math.
141 (2005), no. 1, 146--160.

[K] S. Kleiman, Problem 15. Rigorous fundation of the Schubert's enumerative
calculus, Proceedings of Symposia in Pure Math., 28 (1976), 445-482.

[KK] B. Kostant and S. Kumar, The nil Hecke ring and the cohomology of $G/P$
for a Kac-Moody group G, Adv. Math. 62(1986), 187-237.

[L$_{1}$] J. Leray, Sur l'anneau d'homologie de l'espace homog\`{e}ne,
quotient d'un groupe clos par un sousgroupe ab\'{e}lien, connexe, maximum, C.
R. Acad. Sci. Paris 223, (1946). 412--415.

[L$_{2}$] J. Leray, Propri\'{e}t\'{e}s de l'anneau d'homologie de la
projection d'un espace fibr\'{e} sur sa base, C. R. Acad. Sci. Paris 223,
(1946). 395--397.

[L] L. Lesieur, Les problemes d'intersections sur une variete de Grassmann, C.
R. Acad. Sci. Paris, 225 (1947), 916-917.

[LS$_{1}$] A. Lascoux and M-P. Sch\"{u}zenberger, Polyn\^{o}mes de Schubert,
C.R. Acad. Sci. Paris 294(1982), 447-450.

[LS$_{2}$] A. Lascoux and M-P. Sch\"{u}zenberger, Schubert polynomials and the
Littlewood-Richardson rule, Lett. Math. Phys. 10 (1985), 111--124.

[LPR] A. Lascoux, P. Pragacz and J. Ratajski, Symplectic Schubert polynomials
\`{a} la polonaise, Adv. Math. 140(1998), 1-43.

[LR] D. E. Littlewood and A. R. Richardson, Group characters and algebra,
Philos. Trans. Roy. Soc. London. 233(1934), 99-141.

[M] I. G. Macdonald, Symmetric functions and Hall polynomials, Oxford
Mathematical Monographs, Oxford University Press, Oxford, second ed., 1995.

[Ma] L. Manivel, Fonctions sym\'{e}triques, polyn\^{o}mies de Schubert et
lieux de d\'{e}g\'{e}n\'{e}rescence, Cours Sp\'{e}cialis\'{e}s, no. 3, Soc.
Math. France, 1998.

[Mo] D. Monk, The geometry of flag manifolds, Proc. London Math. Soc., 9
(1959), pp. 253--286.

[PR$_{1}$] P. Pragacz and J. Ratajski, A Pieri-type formula for $SP(2m)/P$ and
$SO(2m+1)/P$, C. R. Acad. Sci. Paris Ser. I Math. 317 (1993), 1035--1040.

[PR$_{2}$] P. Pragacz and J. Ratajski, A Pieri-type formula for Lagrangian and
odd orthogonal Grassmannians, J. Reine Angew. Math. 476 (1996), 143--189.

[PR$_{3}$] P. Pragacz and J. Ratajski, A Pieri-type theorem for even
orthogonal Grassmannians, Max-Planck Institut preprint, 1996.

[RS] J. Remmel and M. Shimozono, A simple proof of the Littlewood-Richardson
rule and applications, Selected papers in honor of Adriano Garsia (Taormina,
1994), Discrete Math. 193 (no. 1-3)(1998), 257--266.

[S$_{1}$] F. Sottile, Pieri's formula for flag manifolds and Schubert
polynomials, Ann. Inst. Fourier (Grenoble) 46 (1) (1996), 89--110.

[S$_{2}$] F. Sottile, Four entries for Kluwer encyclopaedia of Mathematics,
arXiv: Math. AG/0102047.

[St] R. Stanley, Some combinatorial aspects of the Schubert calculus,
Combinatoire et repr\'{e}sentation du groupe sym\'{e}trique, Strasbourg
(1976), 217-251.

[Ste] J. Stembridge, Computational aspects of root systems, Coxeter groups,
and Weyl characters, MSJ. Mem. Vol. 11(2001), 1-38.

[TW] H. Toda and T. Watanabe, The integral cohomology ring of $F_{4}/T$ and
$E_{6}/T$, J. Math. Kyoto Univ., 14-2(1974), 257-286.

[W] T. Watanabe, The integral cohomology ring of the symmetric space EVII, J.
Math. Kyoto Univ., 15-2(1975), 363-385.

[Wi] R. Winkel, On the multiplication of Schubert polynomials, Adv. in Appl.
Math. 20 (1)(1998), 73--97.
\end{document}